\documentclass[11pt]{article}
\usepackage[english]{babel}
\usepackage{amsmath}
\usepackage{amsthm}
\usepackage{amsfonts,amssymb,latexsym,graphicx,psfrag}
\usepackage{version}

\newtheorem{thm}{Theorem}[section]
\newtheorem{cor}[thm]{Corollary}
\newtheorem{lem}[thm]{Lemma}
\newtheorem{prop}[thm]{Proposition}

\theoremstyle{definition}

\theoremstyle{remark}

\numberwithin{equation}{section}

\def \x{\times}

\def \m{\mathbb}

\def \<{\langle}
\def \>{\rangle}

\newcommand{\R}{\mathbb{R}}
\newcommand{\N}{\mathbb{N}}

\newcommand{\B}{\mathbb{B}}

\begin{document}
\title{Regularity of the Hamiltonian along Optimal Trajectories}
\author{M. Palladino and R. B. Vinter\thanks{M. Palladino and R. B. Vinter are with Department of Electrical and Electronic Engineering, Imperial College London, Exhibition Road, London SW7 2BT, UK, {\tt\small m.paladino@imperial.ac.uk,  r.vinter@imperial.ac.uk}}}
\maketitle

\begin{abstract}

\noindent
This paper concerns state constrained optimal control problems, in which the dynamic constraint takes the form of a differential inclusion. If the differential inclusion does not depend on time, then the Hamiltonian, evaluated along the optimal state trajectory and the co-state trajectory, is independent of time. If the differential inclusion is Lipschitz continuous, then the Hamitonian, evaluated along the optimal state trajectory and the co-state trajectory, is Lipschitz continuous. These two well-known results are examples of the following principle: the 
Hamiltonian, evaluated along the optimal state trajectory and the co-state trajectory, inherits the regularity properties of the differential inclusion, regarding its time dependence. We show that this principle also applies to another kind of regularity: if the differential inclusion has bounded variation with respect to time, then the Hamiltonian, evaluated along the optimal state trajectory and the co-state trajectory,  has bounded variation. Two applications of these newly found properties are demonstrated. One is to derive improved conditions which guarantee the nondegeneracy of necessary conditions of optimalty, in the form of a Hamiltonian inclusion. The other  application is to derive new, less restrictive, conditions under which minimizers in the calculus of variations  have bounded slope. The analysis is based on a  new, local, concept of differential inclusions that have bounded variation with respect to the time variable, in which conditions are imposed on the multifunction involved, only in a neighborhood of a given state trajectory.
\end{abstract}
{\small 

\noindent

}


\newpage
\section{Introduction}

Consider the optimal control problem
\ \\

\noindent
$(P)\:\left\{ \begin{array}{l}
\mathrm{Minimize }\, g(x(S),x(T))+\int_{S}^{T}L(t,x(t),\dot{x}(t))dt\\
\mbox{over absolutely continuous functions } x(.):[S,T]\rightarrow \mathbb{R}^{n} \mbox{ s.t. }
\\ 
\dot{x}(t)\,\in \,  F(t,x(t))\;\mathrm{a.e.,}
\\
x(t) \in A
\quad \mbox{for all } t \in [S,T]
\\
(x(S), x(T))\in C\;,
\end{array}\right.
$
\ \\

\noindent
the data for which comprises  an interval $[S,T]$ $(T >S)$, functions $g:\mathbb{R}^{n}\times\mathbb{R}^{n}\rightarrow\mathbb{R}$,  $L:[S,T]\times \R^{n}\times \R^{n} \rightarrow \R$, 
 closed sets $A \subset \R^{n}$ and $C\subset \R^{n} \times \R^{n}$
and a multifunction $F(.,.):[S,T]\times \R^{n} \leadsto \mathbb{R}^{n}$. 
Notice the presence of the pathwise state constraint, $x(t) \in A$. The state constraint set $A$ is assumed to have the  inequality functional representation:
$$
A\,=\,\{x\,|\,h(x)\,\leq\, 0 \} 
$$
for some function $h(.):\R^{n} \rightarrow \mathbb{R}$. By  allowing $h(.)$ to be a nonsmooth function, we capture within  this formulation a wide range of possible state constraints descriptions. An arbitrary closed state constraint set $ A$  is covered by this formulation, since we can always take  $h(x)= d_{A}(x)$, where $d_{A}(.)$ denotes the distance function form the set $A$ (defined below). But it is convenient to derive necessary conditions directly in terms of $h(.)$, because practical constraints frequently come in functional inequality form, and it is desirable in such cases to develop analytical tools for analysing the optimal control problem expressed directly in terms of the  functions arising from the problem formulation.
\ \\

\noindent
A state trajectory $x(.)$  is an absolutely continuous function that satisfies $\dot{x}(t) \in F(t,x(t))$, a.e.
The state trajectory $x(.)$ is said to be feasible if 
$(x(S),x(T))\in C$ and $x(t) \in A$ for all $t \in [S,T]$.
\ \\

\noindent
We say that the state trajectory  $\bar{x}(.)$ is a {\it minimizer} if it achieves the minimum of $g(x(S),x(T))$ over all feasible state trajectories $x(.)$. It is  called a {\it $L^{\infty}$-local minimizer} if, for some $\delta >0$, 
$$
g(x(S),x(T))\;\geq\; g(\bar{x}(S),\bar{x}(T))
$$
for all feasible state trajectories $x(.)$ such that 
\begin{equation}
\label{local}
||x(.)-\bar{x}(.)||_{L^{\infty}} \leq \delta\,.
\end{equation}
For simplicity in this introduction, though not in the analysis to follow, we assume that $L(.,.,.)\equiv 0$ (no integral cost term).  
\ \\

\noindent
A variety of sets of necessary conditions of optimality are known, under hypotheses which impose Lipschtiz continuity conditions on the data, regarding its $x$-dependence, but which require the data to be merely measurable with respect to the $t$ variable. Typically, these assert the existence of a co-state arc $q(.)$ satisfying conditions generalizing the Euler-Lagrange equation or Hamilton's system of equations, the Weierstrass
condition and the transversality condition. Relevant references include \cite{C76},\cite {I97},\cite {LR94},\cite {M95} (no state constraints), and \cite{LR96} ,\cite{vinterz} (when state constraints are present). We refer also to the monographs \cite{clarke1}, \cite{vinter} for expository treatments of this material.
\ \\

\noindent
If we additionally hypothesize regularity of $F(t,x)$ with respect to the $t$ variable, then it is possible to extract additional information about optimal trajectories, expressed in terms of the Hamiltonian function:
$$
H(t,x,q)\,=\, \sup \{q\cdot v\,|\, v \in F(t,x)  \}\,.
$$
Write $H[.]:[S,T] \rightarrow \R^{n}$ for the Hamiltonian evaluated along $\bar{x}(.)$ and $p(.)$:
$$
H[t]\,=\, H(t, \bar x(t), p(t)) \quad t \in [S,T]\,.
$$
A property of  this nature is  \cite{clarke1}, \cite{vinter} :
\ \\

\noindent
{(I): \it `$t\rightarrow F(.,x)$ is constant (for each $x$)' $\implies$ `$H[.]$ is constant'}
\ \\

\noindent
This condition, which has as precursor the 2nd Erdmann condition in the Calculus of Variations, is referred to as the `constancy of the Hamiltonian' condition for autonomous problems.
\ \\

\noindent 
Now suppose that $F(.,x)$  is Lipschitz continuous. Then the optimal control problem can be reformulated as an autonomous problem in which time is interpreted as an additional state variable. Property (I), applied to the reformulated problem, translates into the following information concerning the original problem: 
\ \\

\noindent
{(II): \it `$\,t\rightarrow F(.,x)$ is Lipschitz continuous (for each $x$)' $\implies$ `$H[.]$ is Lipschitz continuous'}
\ \\

\noindent
The latter property is unremarkable, when the state constraint is absent (i.e. $A =\R^{n}$) since, in this case, $q(.)$ is Lischitz continuous, and the Lipschitz continuity of $t \rightarrow H[.]$ can be deduced directly from the Lipschitz continuity of $F(.,x)$. But the Lipschitz continuity of $H[.]$ is perhaps an unexpected property, when the state constraint is present, since, in this setting,  the costate arc $q(.)$ can be discontinuous.  
\ \\

\noindent 
This brings us to the central question addressed in this paper. Can properties (I) and (II) can be interpreted as part of a  general principle (Q)?\ \\

\noindent
{\it (Q): `\,$H[.]$ inherits the regularity of  $F(.,x)$'.}
\ \\

\noindent
The main contribution of this paper is to extend this principle to a larger regularity class. Specifically
\begin{itemize}
\item[(i):] We show that if $F(.,x),$ has bounded variation, then $H[.]$ also has bounded variation.
\item[(ii):] We relate the cummulative variation function of the $H[.]$ to that of the data.
\item[(iii):] We make two applications of these regularity properties of the Hamiltonian, first, to derive new, less restrictive conditions under which necessary conditions of optimality, in the form of Clarke's Hamiltonian inclusion, are nondegenerate and, second, establish boundedness of the derivatives of optimal state trajectories, under hypotheses that are less restrictive than those earlier imposed.
\end{itemize}
As a first step we need to make clear the meaning of the statement `$t \rightarrow F(t,x)$ has bounded variation'. An obvious approach would be to require
$$
\sup \{  \sup_{x \in X_{0}} \sum_{i=0}^{N-1} d_{H}(F(t_{i+1},x),F(t_{i},x) )\}\,<\, \infty\,,
$$
where the outer supremum is taken over all partitions $\{t_{0}=S, t_{1},\ldots, t_{N}=T\}$ of $[S,T]$. (In this relation, $d_{H}(.,.)$ denotes the Hausdorff distance between sets, and $X_{0}$ is some suitably large ball in $\R^{n}$ containing the values of all state trajectories of interest.) \ \\

\noindent
The concept of a multifunction $C(.):[S,T] \leadsto \R^{n}$ of a scalar variable having bounded variation has previously been encountered in connection with Moreau's `sweeping processes' \cite{Moreau}. A sweeping process is a solution to the differential inclusion
$$
-dx/dt \in N_{C(t)}(x(t)),
$$
in which  $C(.)$ is a given multifunction of a scalar variable, taking values, closed, convex sets. Here $N_{C}(x)$  denotes the normal cone from convex analysis. Existence of solutions, in a weak sense, has been established under the hypothesis that $C(.)$ has finite retraction, which is the bounded variation property of this paper (for multifunctions which do not depend on $x$), when the one-sided distance function (or `retraction') is used to define variation, in place of the  Hausdorff distance function.
\ \\

\noindent
There is additional novelty in this paper, regarding our methodology for deriving regularity properties of the Hamiltonian, which involves applying necessary conditions for multi-stage optimal control problems for differential inclusions, similar to those first derived by Clarke and Vinter \cite{cvmulti1}, \cite{cvmulti2}. This is the first instance, to our knowledge, where multi-stage necessary conditions have been used to derive regularity properties of $H[.]$.
As earlier mentioned, the demonstration that $H[.]$ is Lipschitz continuous when the data is Lipschitz continuous with respect to $t$ is a straightforward extension of the constancy property of $H[.]$ for autonomous problems. Demonstrating that $H[.]$ has bounded variation, and obtaining estimates on the cummulative variation function as is required in some applications, is a much more challenging task. The key idea is to approximate the original  optimal control problem with `bounded variation' data, by a multistage autonomous problem (apart from a small perturbation term contributing Lipschitz time dependence). Necessary conditions are derived for the approximating problem, which take account of its autonomous nature, and the desired regularity properties of the Hamiltonian are obtained by passage to the limit.  The precise formulation of the approximating problem, and also the convergence analysis, make use of techniques first used by Arutyunov et al \cite{AA97}, \cite{AAB94}, for the derivation of refined necessary conditions of optimality which provide information about minimizers in some circumstances when standard necessary conditions are degenerate. 
\ \\

\noindent
{\it Notation:} 
For vectors $x\in \m R^{n}$, $|x|$ denotes the Euclidean length. $B$ denotes the 
closed
unit ball in $\m R^{n}$. 
Given a multifunction $\Gamma(.): \m R^{n} \leadsto \m R^{k}$, the graph of $\Gamma(.)$, written Gr$\,\{\Gamma(.)\}$, is the set $\{(x,v)\in \m R^n\x\m R^k\,|\,v\in \Gamma(x)\}$. Given a set $A \subset \m R^{n}$ and a point $x \in \m R^{n}$, we denote by $d_{A}(x)$ the Euclidean distance of a point $x\in \m R^{n}$ from $A$:
$$
d_{A}(x)\;:=\; \inf\{|x-y| \,|\, y \in A \}\;.
$$ 
Given an interval $I$,  we write $\chi_{I}(t)$ for the indicator function of $I$, which takes values $1$ and $0$ when $t\in I$ and $t \notin I$, respectively. 
\ \\

\noindent
A function $r:[S,T] \rightarrow \R$ of bounded variation on the interval $[S,T]$ has a left limit, written $r(t^{-})$, at every point $t \in (S,T]$ and  a right limit, written $r(t^{+})$, at every point $t \in [S,T)$. We say $r(.)$ is normalized if it is right continuous on $(S,T)$. 
\ \\

\noindent
We denote by $NBV^{+}[S,T]$ the space of increasing, real-valued, normalized functions $\mu(.)$ on $[S,T]$  of bounded variation, vanishing at the point $S$. The total variation of a function $\mu(.)\in NBV^{+}[S,T]$ is written $||\mu||_{\mbox{TV}}$. As is well known, each point $\mu(.)\in NBV^{+}[S,T]$ defines a Borel measure on $[S,T]$. This associated measure is also denoted $\mu$.
\ \\

\noindent
We shall use several constructs of nonsmooth analysis. 
Given a closed set $D\subset \R^{k}$ and a point $\bar{x} \in D$. The {\it limiting normal cone} $N_D(\bar{x})$ of $D$ at $\bar{x}$ is defined to be
\begin{eqnarray*}
N_D(\bar{x}):= \Big\{ \;
p \; | \; \exists \; x_{i} \stackrel{D}{\longrightarrow}\bar{x}
, \; p_i \longrightarrow p \;\; \mbox{  s.t. } \; 
\limsup_{x \stackrel{D}{\rightarrow} x_{i}} \;\frac{p_i \cdot (x - x_i)}{|x - x_i|} \leq 0  \quad \mbox{for all } i \in \N \Big\}
\; .
\end{eqnarray*}
Here, the notation $y' \stackrel{D}{\rightarrow} y$ is employed to indicate that all elements in the convergent sequence $\{y'\}$ lie in $D$.
\ \\

\noindent
Take a lower semicontinuous function $f:\m R^k\to \m R\cup\{+\infty\}$ and a point $\bar x\in \mbox{dom }f \,:=\, \{ x\in \R^k\,|\ f(x) < +\infty\}$  The {\it limiting subdifferential of $f$ at $\bar x$} (termed the {\it subdifferential} in \cite{rockafellar}, \cite{vinter}) is denoted $\partial f(\bar x)$:
\begin{eqnarray*}
\lefteqn{\partial f(\bar{x}) \;:=\;  \Big\{ \xi\,|\, \exists\; \xi_{i} \rightarrow \xi \mbox{ and }
x_{i} \stackrel{\mbox{{\tiny {\rm dom}}}\,f}{\longrightarrow} \bar{x} \mbox{ such that}
}
\\[3mm]
&& \limsup_{x \rightarrow  x_{i}} \frac{ \xi_{i} \cdot (x-x_{i})- \varphi(x)+\varphi(x_{i})}{|x-x_{i}|}
\, \leq  0 \, \mbox{ for all } i \in \N \Big\} \;. 
\end{eqnarray*}
For further information, we refer to the monographs \cite{clarke}, \cite{rockafellar} and \cite{vinter}.

\section{Multifunctions of Bounded Variation}

%
%
%
Take a bounded interval $[S,T]$, a compact set $A\subset \R^{k},$ a multifunction $F(.,.,.): [S,T]\times \R^{n} \times A\leadsto \R^{n}$ and a continuous function $\bar{x}(.): [S,T]\rightarrow \R^{n}$. Generic elements in the domain of $F(.,.,.)$ are denoted by $(t,x,a)$.  
\ \\

\noindent
In this section we define a concept that makes precise the statement `$F(t,x,a)$ has bounded variation with respect to the $t$ variable, along  $\bar{x}(.)$, uniformly with respect to $a\in A$'. If $F(t,x,a)$ is independent of $(x,a)$ and single valued, i.e. $F(t,x,a)=\{f(t)\}$ for some function $f(.):[S,T]\rightarrow \R^{n}$, this concept reduces to the standard notion  `$f(.)$ has bounded variation'.
\ \\

\noindent
For any $t \in [S,T]$, $\delta >0$ and partition 
$
{\cal T}\,=\, \{t_{0}=S,t_{1},\ldots,t_{N-1}, t_{N}=t\}
$
of $[S,t]$, 
define $I^{\delta}({\cal T})\in \R^{+} \cup \{+ \infty\}$ to be 
$$
I^{\delta}({\cal T})\,:=\, 
\sum_{i=0}^{N-1}\sup \left\{ d_{H}(F(t_{i+1},x,a) ,F(t_{i},x,a) )\;|\; x \in \bar x([t_{i},t_{i+1}])+\delta B, a\in A  \right\}\;.
$$
Here, $\bar x([t_{i},t_{i+1}])$ denotes the set $\{\bar x(t)\,|\, t \in  [t_{i},t_{i+1}]\}$.
\ \\

\noindent
Take any $\epsilon >0$.  Let $\eta_{\epsilon}^{\delta}(.): [S,T] \rightarrow \R^{+} \cup \{+ \infty\}$ be the function defined as follows: $\eta_{\epsilon}^{\delta}(S)=0$ and, for $t \in (S,T]$,
$$
\eta_{\epsilon}^{\delta}(t)\;=\; \sup \left\{ I^{\delta}({\cal T})\,|\, {\cal T} \mbox{ is a partition of } [S,t] \mbox{ s.t. } \mbox{diam}\{{\cal T}\}\leq \epsilon 
 \right\}\;.
$$
in which
$$
\mbox{diam} \{ {\cal T} \}\,:=\, \sup \{ t_{i+1}-t_{i} \,|\, i=0,\ldots,N-1 \}\,.
$$
Now define the functions $\eta^{\delta}(.), \eta(.): [S,T] \rightarrow \R^{+} \cup \{+ \infty\}$ to be
\begin{eqnarray}
\label{1.1z}
&&\eta^{\delta}(t)\,:=\, \lim_{\epsilon \downarrow 0} \eta^{\delta}_{\epsilon}(t) \;\mbox{ for } t \in [S,T]\,
\\
\label{1.2z}
&&\eta(t)\,:=\, \lim_{\delta \downarrow 0} \eta^{\delta}(t) \;\mbox{ for } t \in [S,T]\,.
\end{eqnarray}

\noindent
{\bf Definition.} Take a  set $A \subset \R^{k}$, a multifunction $F(.,.,.): [S,T]\times \R^{n} \times A\leadsto \R^{n}$ and a  function $\bar{x}(.): [S,T]\rightarrow \R^{n}$.
We say that {\it $t \rightarrow F(t,.,.)$ has bounded variation along $\bar{x}(.)$} over  $A$, if the function $\eta(.)$ given by (\ref{1.2z}) satisfies
$$
\eta(T)\,< \, +\infty\;.
$$
If $t \rightarrow F(t,.,.)$ has bounded variation along $\bar{x}(.)$ uniformly over  $A$, we refer to the function $\eta(.)$ as  the {\it cummulative variation function} of $t \rightarrow F(t,.,.)$ along $\bar{x}(.)$, uniformly over $A$. We also refer to $\eta^{\delta}_{\epsilon}(.)$ and $\eta^{\delta}(.)$ as the $\delta$-perturbed cummulative variation function and $(\delta,\epsilon)$-perturbed cummulative variation function respectively. 
\ \\

\noindent
If $F(t,x,a)$ does not depend on $a$, we omit mention of the qualifier `uniformly over $A$'. A function $t\rightarrow L(t,.,.)$ is said to have bounded variation along $\bar{x}(.)$ uniformly over  $A$, if the associated multifunction $t \rightarrow \{L(t,.,.)\}$ has this property.
 \ \\

\noindent
It is clear that, for any $t \in [S,T]$, $\delta >0$, $\delta' >0$, $\epsilon >0$, $\epsilon' >0$,  
$$
\delta' \leq \delta \mbox{ and }  \epsilon' \leq \epsilon \implies 
0\leq \eta_{\epsilon'}^{\delta'}(t) \leq  \eta_{\epsilon}^{\delta}(t)\;.
$$
 (This relation is valid even when $\eta_{\epsilon}^{\delta}(t)=+  \infty$, according to the rule  $+\infty \leq +\infty$.)

\noindent
It follows that, for fixed $\delta>0$ and $t \in [S,T]$, the functions $\epsilon \rightarrow \eta^{\delta}_{\epsilon}(t)$ and $\delta' \rightarrow \eta^{\delta'}(t)$ are monotone increasing and bounded below on $(0, \infty)$. The limits appearing in the definitions of $\eta^{\delta}(.)$ and $\eta(.)$ (see (\ref{1.1z}) and (\ref{1.2z})) are therefore well-defined. 
\ \\

\noindent
Assume that $t \rightarrow F(t,.,.)$ has bounded variation along $\bar x(.)$, uniformly over $A$.   Then there exist $\bar{ \delta}>0$ and $\bar{\epsilon}>0$ for which $\eta^{\bar \delta}_{\bar \epsilon}(T) < +\infty$. We list further elementary properites of the accumulative variation functions: for any $\delta \in (0,\bar \delta]$ and $\epsilon \in (0,\bar \epsilon]$,

\begin{itemize}
\item[(a):] $t \rightarrow \eta^{\delta}_{\epsilon }(t)$, $t \rightarrow \eta^{\delta}(t)$ and $t \rightarrow \eta(t)$ are  increasing, finite valued functions, and 
\item[(b):] given any $[s,t]\subset [S,T]$ such that $t-s \leq \epsilon$,
$$
d_{H}(F(t,y,a),F(s,y,a)) \,\leq \,
\eta^{\delta}_{\epsilon}(t)-\eta^{\delta}_{\epsilon}(s), 
$$
$\mbox{ for all }y \in \bar{x}(t')+ \delta \B, t'\in [s,t]\,\mbox{ and } \, a\in A.$
\end{itemize}
%
%
As is well known, an $\R^{n}$-valued function of bounded variation on a finite interval may be discontinuous, but it has everywhere left and right limits and it has at most a countable number points of discontinuity. Any multifunction having bounded variation uniformly along  a given continuous trajectory has similar properties, as described in the following proposition, a proof for whcih appears in the Appendix.
\begin{prop}\label{prop1}
Take  a compact set  $A \subset \R^{k}$, a continuous function $\bar{x}(.):[S,T]\rightarrow \R^{n}$ and a  multifunction $F(.,.,.): [S,T]\times \R^{n}\times A \rightarrow \R^{n}$ which has bounded variation along  $\bar{x}(.)$ uniformly over $A$, and take some $\bar \delta >0$ such that $\eta^{\bar \delta}(T)< + \infty$. Assume that 
\begin{itemize}
\item[(C1)]
$F(.,.,.)$ takes values closed, non-empty sets, $F(.,x,a)$ is measurable for each $(x,a) \in \R^{n}\times A$ and there exists $c>0$ 
such that
\begin{equation}
\label{limsup}
F(t,x,a)\,\subset\,c\B \mbox{ for all $x \in \bar{x}(t)+  \bar{\delta}\B$, $t \in [S,T]$, $a\in A$. }
\end{equation}
\item[(C2)]
There exists a modulus of continuity $\gamma(.): \R^{+} \rightarrow \R^{+}$ such that
\begin{equation}
\label{lcontinuous}
F(t,x,a)\,\subset\, F(t,x',a') + \gamma( |x-x'|+|a-a'|) \B
\end{equation}
\mbox{ for all $x,x' \in \bar{x}(t)+  \bar{\delta}\B$, $t \in [S,T]$ and $a,a'\in A$. }
\end{itemize}
 Take $\delta \in (0,\bar \delta)$. Then, 
\begin{itemize}
\item[(a):]For any $\bar s \in [S,T) $ and $\bar t \in (S,T] $ , the one-sided limits
$$
F(\bar{s}^{+},x,a)\,:=\,\underset{s\downarrow \bar s}{\lim} \,F(s,x,a)
$$ 
and 
$$
F(\bar{t}^{-},y,a)\,:=\,\underset{t\uparrow \bar t}{\lim}\, F(t,y,a)
$$ 
exist for every $x \in \bar{x}(\bar{s})+ \delta \B$, $y \in \bar{x}(\bar{t})+ \delta \B$ and $a\in A$.  
\item[(b):]
For any $\bar s \in [S,T) $ and $\bar t \in (S,T] $ 
$$
\underset{s\downarrow \bar s}{\lim} 
\underset{\underset{a\in A} {x \in \bar{x}(\bar s)+\delta B,}}{\sup}
d_{H}(F(\bar{s}^{+},x,a), F(s,x,a))\,=\,0
$$ 
and 
$$
\underset{t\uparrow \bar t}{\lim} 
\underset{\underset{a\in A}{x \in \bar{x}(\bar t)+\delta B,}}{\sup}
d_{H}(F(\bar{t}^{-},x,a), F(t,x,a))\,=\,0
$$ 
\item[(c):]
There exists a countable set ${\cal A}$ such that, for every $t \in (S,T)\backslash {\cal A}$ and $x \in \bar{x}(t)+ \bar \delta \B$,
$$
\underset{t' \rightarrow t}{\lim} \;\,
\underset{\underset{a\in A}{x \in \bar{x}(t)+\bar \delta B}}{\sup}
d_{H}(F(t',x,a), F(t,x,a))\,=\,0\,.
$$ 

\end{itemize}
\end{prop}

\noindent
The next proposition, a proof for which appears in the Appendix,  relates the cummulative variation function of the multifunction $F(.,x,a)$ to that of the derived multifunction $\tilde F(.,x,a)$, obtained by replacing the end-point values by their left and right limits.
\begin{prop}
\label{prop2.2}
Take  a compact set  $A \subset \R^{k}$, a continuous function $\bar{x}(.):[S,T]\rightarrow \R^{n}$ and a  multifunction $F(.,.,.): [S,T]\times \R^{n}\times A \rightarrow \R^{n}$ which has bounded variation along  $\bar{x}(.)$ uniformly over $A$.  Denote by $\eta(.)$  its cummulative variation function. Assume that hypotheses (C1) and (C2) of Prop. \ref{prop1} are satisfied for some $\bar \delta >0$ such that $\eta^{\bar {\delta}}(T) < +\infty$. Take $\delta \in (0, \bar \delta)$ and let $\tilde{F}(.,.,.): [S,T]\times \R^{n}\times A\rightarrow \R^{n}$ be a multifunction such that, for $(t,x,a)\in [S,T] \times \R^{n}\times A$,
\begin{equation}
\label{modified}
\tilde{F}(t,x,a)\,=\,
\left\{
\begin{array}{ll}
F(S^{+},x,a)& \mbox{if } t=S \mbox{ and } |x-\bar{x}(S)| \leq \delta
\\
F(T^{-},x,a)& \mbox{if } t=T \mbox{ and }  |x-\bar{x}(T)| \leq \delta
\\
F(t,x,a)& \mbox{otherwise . } 
\end{array}
\right.
\end{equation}
(Note that limit sets $F(S^{+},x,a)$ and $F(T^{-},x,a))$ exist, by the preceding proposition.) Then  $\tilde{F}(.,.,.)$ has bounded variation along $\bar{x}(.)$ and the 
cummulative variation function $\tilde \eta(.)$ is left continuous at $S$ and right continuous at $T$ respectively, i.e.
$$
\tilde \eta(S)=\lim_{s\downarrow S}\tilde \eta(s)\quad \mbox{and} \quad \tilde \eta(T)=\lim_{t \uparrow T} \tilde\eta(t)\,.
$$
The cummulative variation function of $F(.,.,.)$ and $\tilde F(.,.,.)$ 
are related as follows:
$$
\tilde \eta(t) -\tilde \eta(s) \,=\,  \eta(t) - \eta(s) \quad \mbox{for } [s,t] \subset (S,T),
$$
Furthermore,
$$
\lim_{s\downarrow S} \eta(s) -
\lim_{s\downarrow S}\tilde \eta(s) =\underset{a\in A} {\sup}\, d_{H}(F(S,\bar{x}(S),a), F(S^{+},\bar{x}(S),a) )
$$
and
$$
\lim_{t \uparrow T} \eta(t) -
\lim_{t \uparrow T}\tilde \eta(t) = \underset{a\in A}{\sup}\, d_{H}(F(T,\bar{x}(T),a), F(T^{-},\bar{x}(T),a))\;.
$$
%
%
%
\end{prop}

\ \\

\noindent
The following proposition, a proof for which appears in the Appendix, provides information about how the cummulative variation function of a multifunction $F(.,x,a)$ is affected when the parameter space for $a$ changes. 
\begin{prop}
\label{simple}
Take compact sets $A_{1}\subset A \subset \R^{k}$, a continuous function $\bar{x}(.)$ and a multifunction $F(.,.,.): [S,T]\times \R^{n}\times \R^{k}\leadsto \R^{n}$. Assume that $t \rightarrow F(t,.,.)$ along $\bar{x}(.)$ uniformly over $A$. Write $\eta_{A}(.)$ and $\eta_{A_{1}}(.)$ for the cummulative variation functions of $t \rightarrow F(t,.)$ with respect to the sets $A$ and $A_{1}$ respectively. (Since $t \rightarrow F(t,.,.)$ has bounded variation with respect to $A$, it automatically has bounded variation with respect to the small set $A_{1}$.). Then, for any  $[s,t] \subset [S,T]$,
\begin{equation}
\label{1star}
\eta_{A_{1}}(t) - \eta_{A_{1}}(s) \leq \eta_{A}(t) - \eta_{A}(s)\;.
\end{equation}
\end{prop}


\section{Main Results} 
We  refer to the following hypotheses on the data for problem $(P)$ of the introduction, in which $\bar{x}(.)$ is the $L^{\infty}$-local minimizer of interest and $\bar \delta >0$ is some constant:
\ \\

\noindent
{\bf (H1):}  $F(.,.)$ takes values closed, non-empty sets. For each $x \in \R^{n}$, $F(.,x)$ is a ${\cal L}$ measurable multifunction and $L(.,x,.)$ is a ${\cal L}\times {\cal B}^{n}$ measurable function, where ${\cal L}$ and  ${\cal B}^{n}$ 
denote the Lebegue subsets of $[S,T]$ and the Borel subsets of $\R^{n}$ respectively. $g(.,.)$ is Lipschitz continuous on $(\bar{x}(S),\bar{x}(T) )+\bar \delta \,(\B \times \B)$. 
\ \\

\noindent
{\bf (H2):} $L(.,.,.)$ is bounded on bounded sets and there exist $k>0$ and $c>0$
such that
\ \\

\noindent
$\quad \quad F(t,x)\subset F(t,x')+k(|x-x'|)B\quad$ and
$ \quad F(t,x)\in c \B $
\ \\

\noindent
$\quad \quad|L(t,x,v)-L(t,x',v)| \,\leq \, k|x-x'|
$
\ \\

\noindent
$\quad$for all $x,x'\in\bar{x}(t)+\bar \delta B,\, v\in F(t,x),\,\mathrm{a.e.}\; t\in[S,T].$
\ \\

\noindent
{\bf (H3):} 
There exists a constant $k_{h}$  such that
$$
|h(x)-h(x')| \leq k_{h}|x-x'|
$$
for all $x,x'\in \bar{x}([S,T])+ \bar \delta B$ .
\ \\

\noindent
{\bf (H4):} $\{ (e,v)\in \R^{n}\times \R\,|\, e \in F(t,x),\, v \geq L(t,x,e) \}$ is a convex set
\ \\

\noindent
for each $x \in \R^{n}\,$.
\ \\

\noindent
Define the Hamiltonian function $H_{\lambda}(.,.,.): [0,1]\times \R^{n}\times \R^{n}\rightarrow \R$
\begin{equation}
\label{hamiltonian}
H_{\lambda}(t,x,p):=\max_{v \in F(t,x)} \left[p \cdot v -\lambda L(t,x,v)\right]
\end{equation}
and write
$$
\tilde{C}\,:=\, C \cap\left( \{x' \,|\,  h(x') \leq 0 \}\times \{x' \,|\,  h(x') \leq 0 \} \right)\,.
$$

\begin{thm}
\label{thm:A}
Let $\bar{x}(.)$
be an $L^\infty$-local minimizer for $(P)$. Assume that hypotheses $(H1)$-$(H4)$ are satisfied. Assume also that
\ \\

\noindent
${\bf (BV)}$:
\vspace{0.05 in}

 $t \rightarrow F(t,.)$ has bounded variation along $\bar{x}(.)$
\vspace{0.05 in}

\noindent
and
\vspace{0.05 in}

$t \rightarrow L(t,. , .)$ has bounded variation along $\bar{x}(.)$ uniformly over $c \B$,
\vspace{0.05 in}

\noindent
where $c$ is the constant of hypothesis (H2). Write $\eta_{F}^{*}(.)$ and $\eta_{L}^{*}(.)$ for the normalized cummulative variation functions  of $F(.,.)$  and $L(.,.,.)$ respectively, along $\bar{x}(.)$. 
\ \\

\noindent
Then 
there exist an absolutely continuous function $p(.):[S,T] \rightarrow  R^n$, a function $\mu(.)\in NBV^{+}[S,T]$, a $\mu$-integrable function $\gamma(.)$ 
and $\lambda \geq 0$ such that
\begin{itemize}
\item[(i):] $(p(.), \lambda, \mu(.))\; \not= \;(0,0,0)$,
\item[(ii):] $(-\dot p(t),\dot{\bar{x}}(t))\in \tilde{\partial}_{x,p} H_{\lambda}(t,\bar{x}(t),q(t)) $ a.e.,
\item[(iii):] $q(t)\cdot \dot{\bar{x}}(t)- \lambda L(t, \bar{x}(t), \dot{\bar{x}}(t))\,=\;$
\item[]
\hspace{0.5 in}$ \max \,\{ q(t)\cdot v- \lambda L(t, \bar{x}(t),v) \,|\, v \in F(t, \bar{x}(t)) \}$ a.e.,
\item[(iv):] 
$(q(S),-q(T)) \in  \lambda \partial g(\bar x(S),\bar x(T))+ N_{\tilde{C}}(\bar x(S),\bar x(T))$,
\item[(v):]
 $\gamma(t)\in  \mbox{co}\, \partial_{x}^{>}h(\bar{x}(t))\quad \mu\mbox{-}\mathrm{a.e.}\; t\in[S,T]\,.$
\end{itemize}

\noindent
where
$$
q(t)\,=\, 
\left\{ 
\begin{array}{ll}
p(t) +  \int_{[S,t)}\gamma(s)\mu(ds) & \mbox{if } 
 t \in [S,T)
\\
p(T)+\int_{[S,T]}\gamma(s)\mu(ds) & \mbox{if } t = T\,.
\end{array}
\right.
$$

\noindent
Furthermore, the Hamiltonian evaluated along $(\bar{x}(.),q(.))$ has the following properties: there exists a normalized function of bounded variation $r(.):[S,T]\rightarrow \R $ which is right and left continuous at $S$ and $T$ respectively, i.e.
$$
r(S)= \lim_{s\downarrow S} r(s)\,,\quad  r(T)= \lim_{t\uparrow T} r(t)\,,
$$
such that
\begin{itemize}
\item[(vi):]
$
|r(t)-r(s)|\,\leq\,  ||q(.)||_{L^{\infty}} \times \left(\eta_{F}^{*}(t)-\eta_{F}^{*}(s) \right)
+ \lambda \times \left( \eta_{L}^{*}(t)-\eta_{L}^{*} (s) \right)
$

\noindent
for all $[s,t]\subset (S,T)$\,,
%
\item[(vii):]
$H_{\lambda}(t,\bar{x}(t),q(t))\,=\, r(t) \quad \mbox{a.e.}
$
and
\item[(viii):]
$
r(S)= \lim_{s\downarrow S} H_{\lambda}(s, \bar{x}(S),p(S))\,,\;  r(T)= \lim_{t\uparrow T}H_ \lambda (t, \bar{x}(T),q(T))\;.
$
\end{itemize}
\end{thm}
\ \\

\noindent
In the theorem statement, the `hybrid' partial subdifferentials $\partial^{>}_{x}h(x)$ and $\tilde{\partial}_{x,p}H(s,x,p)$ are defined to be:
%
$$ \partial^{>}_{x}h(x)\;:=\; \{\xi \,| \, \exists\,\,
x_{i}\rightarrow x,  \xi_{i}\rightarrow \xi, \mbox{ s.t. } \xi_{i} = \nabla_{x}h (x_{i}) \mbox{ and } h(x_{i}) > 0  \mbox{ for each  }i\} \;
$$

\noindent
and
\begin{equation}
\label{hybridgradient}
\tilde{\partial}_{x,p}H(t,x,p)\,:=\, \underset{s \rightarrow t}{\mbox{lim sup}} \;\mbox{co}\, \partial_{x,p} H(s,x, p)\,.
\end{equation}
\noindent
{\bf Comments}
\begin{itemize}
\item[(a):] 
The standard Hamiltonian inclusion  $(-\dot p(t),\dot{\bar{x}}(t))\in \mbox{co}\, \partial_{x,p} H_{\lambda}(t,\bar{x}(t),q(t))$ involving the Clarke subdifferential $ \partial_{x,p} H$ implies the Weierstrass condition (iii).  When the Hamiltonian inclusion takes the weaker form (ii) involving the hybrid subdifferential $ \tilde{\partial}_{x,p} H$, the Weierstrass condition (iii) does not automatically follow and is included in the theorem statement as a separate condition.
\item[(b):] Suppose  $p(t)+ \int_{[S,t]}  \nabla h(\bar{x}(s)))ds$ has bounded variation along $\bar{x}(.)$. Then, under hypotheses (H1)-(H4), it is easy to show directly that, for any $p(.), \bar{x}(.)  \in W^{1,1}$ and measure $\mu(.)$, the function 
$$
H(t):= t \rightarrow \max\{(p(t)+ \int_{[S,t]}  \nabla h(\bar{x}(s)))ds)\cdot v \,| v \in F(t,\bar{x}(t))\,   \}
$$ 
 also has bounded variation. Notice that, in this simple direct demonstration, the cummulative variation of both $p(t)+ \int_{[S,t]}  \nabla h(\bar{x}(s)))ds$ and of $p(t)+ \int_{[S,t]}  \nabla h(\bar{x}(s)))ds$ both contribute to the cummulative variation function of $h(t)$. So the assertion, merely, that the Hamiltonian is of bounded variation is a trivial addition to known necessary conditions.  However the theorem contributes the extra information (vi), when may be paraphrased as the assertion:  when $(p(.), \mu)$ are multipliers associated with a $L^{\infty}$ minimizer $\bar{x}(.)$, then  the cummulative variation function of $h(.)$ is absolutely continuous with respect to the cummulative variation functions of  $t \rightarrow F(t,.)$ and $t \rightarrow L(t,.,.)$. This tells, perhaps surprisingly, that  {\it only the cummulative varation of  $t \rightarrow F(t,.)$ and $t \rightarrow L(t,.,.)$, and not that of $p(t)+ \int_{[S,t]}  \nabla h(\bar{x}(s)))ds$, contribute to the cummulative variation of the Hamiltonian}. This is a highly nontrivial addition to the standard necessary conditions, even when there are no state constraints. We investigate  implications in the following two sections.
\end{itemize}
\section{Application 1: Minimizer Regularity}
Take a function $L:[S,T]\times \R^{n}\times \R^{n} \rightarrow \R$ and points $x_{0},x_{1} \in \R^{n}$. Consider the optimization problem:
\ \\

\noindent
$(Q)\:\left\{ \begin{array}{l}
\mathrm{Minimize }\, \int_{S}^{T}L(t,x(t),\dot{x}(t))dt\\
\mbox{over absolutely continuous functions } x(.):[S,T]\rightarrow \mathbb{R}^{n} \mbox{ s.t. }

\\
x(S)=x_{0} \; \mbox{and},\; x(T)= x_{1}\;.
\end{array}\right.
$
\ \\

\noindent
In this section we write $H(t,x,p)$ in place of $H_{\lambda}(t,x,p)$ when $\lambda=1$, thus 
$$
H(t,x,p)\,:=\, \sup \{ p\cdot v  - L(t,x,v)\,|\, v \in \R^{n}  \}\,.
$$
It is well-known that $(Q)$ has a minimizer $\bar{x}(.)$ under the following hypotheses:
\begin{itemize}
\item[(HE):]
\begin{itemize}
\item[(i):] $L(.,.,.)$ is ${\cal L}\times {\cal B}^{n\times n}$ measurable, where ${\cal L}$ and ${\cal B}^{r}$ denote the Lebesgue subsets of $[S,T]$ and the Borel subsets of $\R^{r}$ respectively, and $L(t,.,.)$ is lower semicontinuous for each $t \in [S,T]$.
\item[(ii):] $L(t,x,.)$ is convex for each $(t,x) \in  \R^{n}$.

\item[(iii):] There exists a convex function $\theta(.):\R^{+} \rightarrow \R^{+}$ and a number $\alpha $ such that
$$
\lim_{r \uparrow \infty}\; \theta(r)/r \,=\, + \infty,\; 
$$
and
$$
L(t,x,v)  \geq \theta(|v|) -\alpha |x| \mbox{ for all }(t,x,v) \in [S,T]\times \R^{n}\times \R^{n}\,.
$$
\end{itemize} 
\end{itemize}
If (HE) is supplemented by the mild regularity hypothesis on $L(t,.,.)$:
\begin{itemize}
\item[(HR):]  Given any $D>0$ there  exists $k_{D} >0 $ such that

 $L(t,x,v)-L(t,x',v')\,\leq\, k_{D}(|(x,v)-(x',v')|)$

$\hspace{1.4 in} \, \mbox{ for all } (x,v),(x',v')\in D\B,\, t \in [S,T]$
\end{itemize}
which permit us to calculate subdifferentials of $L(t,.,.)$ in some sense,
then we might expect that the minimizer $\bar{x}(.)$ satisfies standard first order conditions of optimality. 
This however is not the case. There are counter-examples (\cite{BM85},\cite{CV84}) of problems with smooth data, satisfying (HE) and (HR), whose minimizers fail to satisfy the Euler Lagrange equation. 
\ \\

\noindent
There has been longstanding interest in the question of what {\it additional} hypotheses, besides, (HE) and (HR), are required to ensure minimizers satisfy classical necessary conditions (or their modern nonsmooth analogues). In much of this literature, additional hypotheses are imposed to ensure that minimizers are Lipschitz continuous, since the Hamiltonan inclusion (for example) is automatically satisfied by Lipschitz continuous minimizers under (HE) and (HR).
One such additional hypothesis (see $[ ]$), guaranteeing that minimizers are Lipschitz continuous and therefore satisfy necessary conditions of interest is
\begin{itemize}
\item[(HA):] $L(.,.,.)$ is locally Lipschitz continuous (in all variables), $L(.,.,.)$ continuously differentiable on a neighborhood of $(t,\bar{x}(t),\dot{\bar{x}}(t))$ for each $(t,x,v) \in (S,T) \times\R^{n}\times \R^{n}$, and 
$$ 
t \rightarrow \nabla_{t}L(t, \bar{x}(t), \dot{\bar{x}}(t) )\mbox{ is integrable}.
$$
\end{itemize}
(For simplicity, we do not state the more complicated form this condition takes, when $L(.,.,.)$ is no longer continuously differentiable.) 
Hypothesis (HA), or to be precise the version of it which does not require continuous differentiability, is automatically satisfied  when $L(t,x,v)$ does not depend on $t$.
\ \\

\noindent
The theory of previous sections will permit us to replace the hypothesis (HA) by a weaker hypothesis requiring that $t \rightarrow L(t,.,.)$ is merely {\it of bounded variation}  is some uniform sense. This new hypothesis covers different cases and it allows $L(t,.,.)$ to have a countable number of fractional singularities, or even discontinuities. 
Specifically, we shall assume: 
\begin{itemize}
\item[(BV):]  There exists $\epsilon >0$ and $\delta >0$  and $K >0$ such that, for all
 partitions 
$
{\cal T}\,=\, \{t_{0}=S,t_{1},\ldots,t_{N-1}, t_{N}=T\}
$
of $[S,T]$  with  diam$\,\{{\cal T}\} \leq \epsilon$, we have
$$
\sum_{i=0}^{N-1}\sup \left\{ |L(t_{i+1},x,v) - L(t_{i},x,v)| \;|\; x \in \bar x([t_{i},t_{i+1}])+\delta B,v \in \R^{n}  \right\} \;\leq\; K\;.
$$
\end{itemize}
(In the language of previous sections, (BV) is the requirement that $t \rightarrow L(t,.,.)$ has bounded variation along $\bar{x}(.)$ uniformly over the set $A$, where now $A$ is the {\it unbounded} set $\R^{n}$.)
\begin{prop}
\label{prop:A1}
Let $\bar{x}(.)$
be an $L^\infty$-local minimizer for $(Q)$. Assume that hypotheses $(HE)$ and $(BV)$ are satisfied.
Then $\bar{x}(.)$ is Lipschitz continuous and there exists a Lipschitz continuous function $p(.): [S,T] \rightarrow \R^{n}$ such that
\begin{equation}
\label{cv}
(-\dot p(t),\dot{\bar{x}}(t))\in \tilde{\partial}_{x,p} H(t,\bar{x}(t),p(t))  \mbox{ a.e. } 
\end{equation}
(where the subdifferential $\tilde{\partial}_{x,p} H(.,.,.)$ is as in (\ref{hybridgradient})).
\end{prop}

\noindent
{\bf Proof.} Fix any $k' > ||\dot{\bar{x}}(.)||_{L^{\infty}
([\bar{t},t'];\R^{n})}$.  Then, under the hypothesis (BV), the function $t \rightarrow L(t,.,.) $ has bounded variation along $\bar{x}$ uniformly over the compact set $k' \B$. Write the cummulative variation function $t \rightarrow \eta(t;k')$. It is clear from the hypothesis that 
$$
\eta(t; k')\,\leq \, K\quad \mbox{for all } t \in [S,T] \mbox{ and }k' \geq 0\;.
$$
This inequality, in which $K$ is the constant of hypothesis (BV), continues to be satisfied if $\eta(t; k')$ is replaced by its normalized version $\eta^{*}(t; k')$, obtained by replacing values of the function at interior points by their right limits.
\ \\

\noindent
We deduce from Tonelli regularity theory (see, e.g.,  \cite[Chapter XX]{vinter}) that there exists a point $\bar{t} \in (S,T)$ such that $\dot{\bar{x}}(.)$ is essentially bounded on a neighborhood of $\bar{t}$. We shall show that $\dot{\bar{x}}(.)$ is essentially bounded on $[\bar t, T]$. An analogous argument `in reverse time' will tell us that  $\dot{\bar{x}}(.)$ is essentially bounded also $[\bar t, T]$. It will follow that $\bar{x}(.)$ is Lipschitz continuous on all of $[S,T]$.
\ \\

\noindent
The proof is based on the fact:
\ \\

\noindent
{\it If there exists   $\bar{k}$ such that, for any $[\bar t, t'] \subset [\bar t, T]$, $t' > \bar{t}$, on which $\dot{\bar{x}}(.)$ is essentially bounded, there exists $p(.)\in W^{1,1}([\bar{t},T];\R^{n})$ such that
\begin{eqnarray}
&&(-\dot{p}(t),\dot{\bar{x}}(t)) \,\in\, \tilde \partial_{x,p} H(t, \bar{x}(t), p(t))\; \mbox{ a.e.} \hspace{2.5 in}
\label{inclusion}
\\
&&|p(t)|\,\leq \,  \bar{k} \; \mbox{ for all } t \in [\bar{t},T]\,.
\nonumber
\end{eqnarray}
Then $\dot{\bar{x}}(.)$ is essentially bounded on $[\bar{t}, T]$.
} 
\ \\

\noindent
(This is a slight modification of \cite[Lemma 11.4.1]{vinter}, in which the right side of (\ref{inclusion}) is expressed in terms of the hybrid subdifferential $\tilde \partial H$; this makes no essential difference to the proof.)
\ \\

\noindent
Take any  $[\bar t, t'] \subset [\bar t, T]$, $t' > \bar{t}$, on which $\dot{\bar{x}}(.)$ is essentially bounded. Now apply Thm. \ref{thm:A} to $\bar{x}(.)$ restricted to $[\bar{t},t']$, regarded as minimizer for the problem
$$
\mbox{Minimize } \left\{ 
\int_{\bar{t}}^{t'} L(s,x(s), \dot{x}(s))ds\,|\, x(\bar{t}) =  \bar x(\bar{t}), x(t')= \bar x(t') \mbox{ and } |\dot x(.)| \leq k'
\right\}\,.
$$
We deduce the existence  of $p(.) \in W^{1,1}$ and  a normalized function of bounded variation (on the interval $[\bar t, t']$) such that
\begin{eqnarray}
&&(-\dot p(t),\dot{\bar{x}}(t))\,\in\, \tilde{\partial}_{x,p} H(t,\bar{x}(t),q(t)) \; \mbox{ a.e.}
\label{6.1}
\\
&& r(t)= \max_{v \in k'\B} \{p(t)\cdot v - L(t,\bar{x}(t), v)\},\; \mbox{a.e.}
\label{6.2}
\\
&&
 p(t) \cdot \dot{\bar{x}}(t) -  L(t,\bar{x}(t), \dot{\bar{x}}(t))\,\geq\, p(t)\cdot v - L(t,\bar{x}(t), v)\; \mbox{ for all } v \in k' \B
\label{6.2aaa}
\\
&& r(t)-r(s)\,\leq \, \eta^{*}(t; k') -\eta^{*}(s;k')\,\mbox{ for all } [s,t]\subset [\bar t, t']\;.
\label{6.2aa}
\end{eqnarray}
Notice the necessary conditions are stated in `normal form', i.e. with cost multiplier $\lambda$ set to $1$. This is permissible since,  for fixed endpoint problems with  a dynamic constraint $\dot{x} \in k' \B$  and minimizing arcs $\bar x(.)$ with velocities interior to $k'\B$, the necessary conditions  above can only be satisfied by multipliers $(p(.),\lambda)$ for which $\lambda >0$. We may then arrange $\lambda=1$ by scaling $p(.)$ and $\lambda$ appropriately.
\ \\

\noindent
Since $\dot{\bar{x}}(t)$ interior to $k' \B$ a.e., and $L(t,\bar{x}(t),.)$ is convex, (\ref{6.2aaa}) implies the `global' optimality property
\begin{equation}
\label{6.3}
p(t) \cdot \dot{\bar{x}}(t)-  L(t,x(t), \dot{\bar x}(t))\geq p(t)\cdot v-   L(t,x(t), v) \; \mbox{ for all $v \in \R^{n}$, a.e. }
\end{equation}

\noindent
Since, additionally, $L(t,\bar{x}(t), .)$ is locally Lipschitz continuous, we have
$$
p(t) \,\in\,   \partial_{v} L(t, \bar{x}(t), v=\dot{\bar{x}}(t)) \; \mbox{ a.e.}
$$
It follows then from hypothesis (HR) that, for some $t_{1} \in (\bar{t},t')$, 
\begin{equation}
\label{6.4}
|p(t)|\,\leq \, k_{1}\quad \mbox{a.e. } t  \in [\bar t, t_{1}]
\end{equation} 
for some $k_{1}>0$ that does not depend on $t'$. We can deduce from (\ref{6.2}) and (\ref{6.4}) that
\begin{equation}
\label{6.5}
|r(s)| \leq k_{2}
\quad \mbox{a.e. } t  \in [\bar t, t_{1}]\,,
\end{equation} 
for some $k_{2}>0$ that does not depend on $t'$. Now write
$$
k_{3}\,:=\,  \underset{t \in [S,T]}{\mbox{max}}\underset{v \in \B}{\mbox{max}} \;L(t, \bar x(t), v)\;,
$$
By considering the choice $v = p(t)/|p(t)|$ if $p(t)\not= 0$ or $v=0$ otherwise in (\ref{6.3}), we deduce that
$$
|p(t)|\,\leq \, k_{3}+ r(t)\quad \mbox{a.e.}
$$
Now choose any $s \in [\bar t, t_{1}]$ such that  (\ref{6.5}) is valid. Then for any $t \in [t_{1},t']$ we have, by (\ref{6.2aa}),
$$
p(t)\,\leq\, k_{3}+r(t)= k_{3}+r(s)+r(t)-r(s)\,\leq\, k_{2}+ k_{3} + \eta^{*}(t;k') -\eta^{*}(s,k') \leq k_{2}+ k_{3}+ K,
$$
where $K$ is the constant in hypothesis (BV). Combining this relation with (\ref{6.4}) and noting that $p(.)$ is continuous, we conclude that the values of  $p(.)$ are uniformly bounded on all of $[\bar{t},t' ]$, and the bound does not depend on $t'$.
 But, as we have earlier observed, the existence of such a bound implies $\dot{\bar{x}}(.)$ is essentially bounded on $[\bar t, t' ]$. 
\ \\

\noindent
We have shown that $\bar{x}(.)$ is Lipschitz continuous. To complete the proof we have merely to recall that, under (HE) and (HR), Lipschitz continuous minimizers satisfy the asserted necessary conditions.

\section{Application 2: Non-Degeneracy of the State Constrained Hamiltonian Inclusion.}
Consider the optimal control problem
\ \\

\noindent
$(P_{1})\:\left\{ \begin{array}{l}
\mathrm{Minimize }\, g(x(T))\\
\mbox{over absolutely continuous functions } x(.):[S,T]\rightarrow \mathbb{R}^{n} \mbox{ s.t. }
\\ 
\dot{x}(t)\,\in \,  F(t,x(t))\;\mathrm{a.e.,}
\\
h(x(t))\,\leq\, 0 \quad \mbox{for all } t \in [S,T]
\\
x(S)=x_{0}, 
\;,
\end{array}\right.
$
\ \\

\noindent
in which $x_{0}$ is a given $\R^{n}$ vector.
To simplify the subsequent analysis, we assume:
\begin{itemize}
\item[$(H3)'$:] $h(.)$ is a continuously differentiable function satisfying $\nabla h(x_{0})\not=0$.
\end{itemize}

\noindent
$(P_{1})$ will be recognized as a special case of $(P)$ in which $C=\{ x_{0}\} \times \R^{n}$ (`fixed initial state'). Let $\bar{x}(.)$ be an $L^{\infty}$-local minimizer for $(P_{1})$.  The state constrained Hamiltonian inclusion condition tells us: under unrestrictive hypotheses (details of which do not concern us here), there exist $p(.)\in W^{1,1}\left([S,T];\m R^n\right)$, a function $\mu(.)\in NBV^{+}[S,T]$
and $\lambda \geq 0$ such that
\begin{itemize}
\item[(i)] $(p(.), \lambda, \mu(.))\; \not= \;(0,0,0)$,
\item[(ii)] $(-\dot p(t),\dot{\bar{x}}(t))\in \tilde{\partial}_{x,p} H(t,\bar{x}(t),q(t)) $ a.e. ,
\item[(iii)] 
$-q(1) \in  \lambda \partial g(\bar x(T))$
\end{itemize}

\noindent
and
\begin{itemize}
\item[(iv)]
 $\mbox{supp}\, \{d\mu\}\subset \{t\,|\, h(t,\bar{x}(t))=0   \} $
\end{itemize}

\noindent
where
$$
q(t)\,=\, 
\left\{ 
\begin{array}{ll}
p(t)+\int_{[S,t)}\nabla h(x(s))\mu(ds) & \mbox{if } t \in [S,T)
\\
p(t)+\int_{[S,T]}\nabla h(x(s))\mu(ds) & \mbox{if } t =T\,.
\end{array}
\right.
$$
These conditions are deficient when 
\begin{equation}
\label{boundary}
h(x_{0})=0 \quad \mbox{`initial state in the boundary'}\,.
\end{equation}
Indeed it is easy to show that, in this case, {\it any} feasible $F$-trajectory $\bar{x}(.)$ satisfies conditions (i)-(iii) above,  when we choose multipliers:
\begin{equation}
\label{degenerate}
p(.)\equiv -\nabla h(\bar{x}(S)),\; \mu(t) = \delta(t) \mbox{ and } \lambda=0\,,
\end{equation}
and so the standard set of conditions (i)-(iii) convey no useful information about minimizers whatsoever.\ \\

\noindent
As the proof of the following proposition shows, the additional information about the Hamiltonian in Thm. \ref{thm:A}, coupled with an inward pointing condition at the left endpoint, provides an extra condition:
$$
\lambda + \int_{(S,T]}d\mu(s)  \,\not=\, 0
$$
This condition excludes the `trivial' set of multiplers (\ref{degenerate}), in which $\lambda=0$ and $\mu(.)$ comprises merely an atom at $\{S\}$, and therefore strengthens the standard set of conditions (i)-(iii). Write
\ \\

\noindent
$$
H(t,x,p)\,=\,\sup\{p\cdot v \,|\, v \in F(t,x) \}\,.
$$
\begin{cor}
Let $\bar{x}(.)$
be an $L^\infty$-local minimizer for $(P_{1})$. Assume that hypotheses $(H1)$, $(H2)$ (for $L(.,.,.)\equiv 0$ ) and (H3$\,'$) are satisfied. Assume also that
\ \\

\noindent
${\bf (BV)}$:
 $t \rightarrow F(t,.)$ has bounded variation along $\bar{x}(.)$
\vspace{0.05 in}

\noindent
and
\vspace{0.05 in}

\noindent
${\bf (I)}$: There exists $v \in \underset{s\downarrow S} {\mbox{lim inf}}\;F(s,x_{0})
$
such that
$$\nabla h(x_{0})\cdot v < 0\;. 
$$
Then 
there exist $p(.)\in W^{1,1}\left([S,T];\m R^n\right)$, a function $\mu(.)\in NBV^{+}[S,T]$
and $\lambda \geq 0$ such that
\begin{itemize}
\item[(i)] $ \lambda + \int_{(S,T]} \mu(ds))\; \not= \;0$,
\item[(ii)] $(-\dot p(t),\dot{\bar{x}}(t))\in \tilde{\partial}_{x,p} H(t,\bar{x}(t),q(t)) $ a.e. ,
\item[(iii)] 
$(q(S),-q(T)) \in  \lambda \partial g(\bar x(S),\bar x(T))+ N_{\tilde{C}}(\bar x(S),\bar x(T))$,
\end{itemize}

\noindent
where
$$
q(t)\,=\, 
\left\{ 
\begin{array}{ll}
p(t)+\int_{[S,t)}\nabla h(x(s))\mu(ds) & \mbox{if } t \in [S,T)
\\
p(t)+\int_{[S,T]}\nabla h(x(s))\mu(ds) & \mbox{if } t =T\,.
\end{array}
\right.
$$
\end{cor}
\vspace{0.1 in}
\noindent
{\bf Proof.} From  Thm. \ref{thm:A} we know that there exist  $p(.)\in W^{1,1}\left([S,T];\m R^n\right)$, a function $\mu(.)\in NBV^{+}[S,T]$
and $\lambda \geq 0$ such that (ii) and (iii) are satisfied, and $||p(.)||_{L^{\infty}} + \lambda+ \int_{[S,T]}d\mu(s) \not= 0$. 
\ \\

\noindent We must check also condition (i). Suppose this condition is not satisfied. Let $\tilde q(.)$ be the modification to $q(.)$, in which $q(S)$ is replaced by its left limit $q(S^{+})$. Then $\tilde q(.)$ is an absolutely continuous funtion satisfying  
$$ (-\dot {\tilde{q}} (t),\dot{\bar{x}}(t)) \in \partial_{x,p} H(t,\bar{x}(t),\tilde q(t))  \quad \mbox{a.e.}
$$
and $\tilde q(T)=0$. It can be deduced, with the help of Gronwall's inequality, that $\tilde q(.) \equiv 0$. So $q(t)=\tilde q(t)=0$ for $t \in (S,T]$. Then
$$
H(t, \bar{x}(t),q(t))=0 \;\mbox{ for all } t\in (S,T] \;.
$$
Now let us examine the implications of the `extra' information about $\bar{x}(.)$ supplied by Thm. \ref{thm:A}: there exists a normalized function of bounded variation $r(.): [S,T]\rightarrow \R$, continuous at its endpoints, such that
$$
r(t)= H(t, \bar{x}(t),q(t))=0 \;\mbox{ a.e } t\in (S,T] \;
\; \mbox{ and }  r(S)= H(S^+,\bar x(S),p(S))
\;.
$$
From the first condition and the continuity of $r(.)$ at S, we conclude that $r(S)=0$. 
This implies that
\begin{equation}
\label{contra}
 \underset{v \in F(S^{+},x_{0})}{\mbox{sup}} p(S)\cdot v \,=\,0\;.
\end{equation}
But $p(S)= \tilde q(S+)- \nabla h(x_{0}) \mu\{\{S\}\} =0-  \nabla h(x_{0})\mu\{\{S\}\}$. There are now two possibilities. It $\mu(.)\equiv 0$ then $p(.)= \tilde q(.)\equiv 0$. This contradicts the non-triviality condition $||p(.)||_{L^{\infty}}+ \lambda + \int_{[S,T]}d \mu(s)) \not=0$. On the other hand, if $\mu\{\{S\}\}\not= \gamma >0$ then 
$$
H(S^+,x_{0},p(S))= \underset{v \in F(S^{+},x_{0})}{\mbox{sup}} p(S)\cdot v= -  \underset{v \in F(S^{+},x_{0})}{\mbox{inf}}\gamma \nabla h(x_0)\cdot v >0\;,
$$
by hypothesis (I). This contradicts (\ref{contra}). This concludes the proof. 
\ \\

\noindent. 
{\bf Comments.}
\begin{itemize}
\item[(a):] The idea of using regularity of the Hamiltonian to justify the exclusion of trivial multipliers in state constraint optimal control is due to Arutyunov and Aseev.  They assume that the data is Lipschitz continuous with respect to time. The corollary improves on this earlier work, by excluding the trivial multiplier set when the data is merely of bounded variation with respect to time.
\item[(b):] The corollary  excludes just one type of non-degeneracy. The condition (i) still permits $\lambda=0$, a case in which the necessary condtions make no reference to the cost function and merely convey information about the constraints. The hypotheses can be strengthened to ensure also normality, i.e. we can choose $\lambda \not= 0$, in cases not covered in the earlier `non-degeneracy' literature. 
\end{itemize}

\section{Preliminary Analysis: Discrete Approximations}
A key step in the proof the $H[.]$ has bounded variation when $F(.,x)$ has bounded variation will be to derive necessary conditions for a perturbation of problem $(P)$, resulting from the discrete approximation of the multifunction $F(.,x)$ (and the addition of penalty terms to the cost). The perturbed problem $(P)'$ has the form of a multistage optimal control problem, with end-times taken to be the grid points of a partition of $[S,T]$: 
\ \\

\noindent
$(P')\:\left\{ \begin{array}{l}
\mathrm{Minimize }\, g(x(S),x(T))+\sum_{i=0}^{N-1} \int_{t_{i}}^{t_{i+1}}L_{i}(x(t),\dot{x}(t))dt +
\int_{S}^{T}e(t,x(t))dt
\\
\mbox{over absolutely continuous functions } x(.):[S,T]\rightarrow \mathbb{R}^{n} \mbox{ s.t. }
\\ 
\dot{x}(t)\,\in \,  F_{i}(x(t))\;\mathrm{a.e. }\; t \in [t_{i},t_{i+1}), \, i=0, \ldots, N-1
\\
(x(S), x(T))\in C\;,
\end{array}\right.
$
\ \\

\noindent
for which the data is as follows:  a partition $\{t_{0}=S,t_{1},\ldots,t_{N-1}, t_{N}=T\}$ of $[S,T]$, multifunctions and functions $F_{i}(.):\R^{n}\leadsto \R^{n}$ and $L_{i}: \mathbb{R}^{n}\times \mathbb{R}^{n}\rightarrow \mathbb{R},\, i=0,\ldots, N-1$,  functions
 $e(.,.):[S,T]\times \R^{n}\rightarrow \R$ and $g(.,.):\R^{n}\times \R^{n}\rightarrow \R$, 
and a set $C\subset \R^{n}\times \R^{n}$. 
The special structure of the cost integrand reflects the requirements of the analysis to follow.   
\ \\

\noindent
Problem $(P)'$ will be recognised as a special case of the problem $(P)$ of the introduction, when  the following identifications are made for $F(.,.)$ and $L(.,.,.)$:
\ \\

$
F(t,x)\,=\,\sum_{i=0}^{N-1}F_{i}(x)\chi_{[t_{i},t_{i+1})}(t)\; \mbox{ a.e.}\,,
$ 
\ \\

$L(t,x,v)\,=\, \sum_{i=0}^{N-1}L_{i}(x,v)\chi_{[t_{i},t_{i+1})}(t)+e(t,x)$\, .
%
\ \\

\noindent
Define the $i$'th stage hamiltonian $H^{i}_{\lambda}(.,.,.)$ (parameterized by the cost multiplier $\lambda$) to be
$$
H^{i}_{\lambda}(t,x,p):= \max_{v \in F(t,x)} \{ p \cdot v - \lambda L(t,x,v) \}
$$
for $(x,p)\in \R^{n}\times \R^{n}$ and $t \in [t_{i}, t_{i+1})$, $i=0,\ldots,N-1$.

\begin{prop}
\label{multistage}
Let $\bar{x}(.)$
be an $L^\infty$-local minimizer for $(P')$. Assume, for some $\bar \delta >0$, 
\ \\

\noindent
{\bf (HM1):}  $F_{i}(.)$ takes values closed, non-empty sets, $i=0,\ldots,N-1$ and $L(.,.)$ is continuous.
 $g(.,.)$ is Lipschitz continuous on $(\bar{x}(S),\bar{x}(T) )+\bar \delta \,(\B \times \B)$. 
\ \\

\noindent
{\bf (HM2):} There exist $k>0$ and $c>0$, such that, for each $i\in \{0,\ldots,N-1\}$, $L_{i}(.,.)$ is bounded on bounded sets, and
\ \\

\noindent
$\quad \quad F_{i}(x)\subset F_{i}(x')+k(|x-x'|)B\quad$ and
$ \quad F_{i}(x)\in c \B $
\ \\

\noindent
$\quad \quad|L_{i}(x,v)-L_{i}(x',v)| \,\leq \, k|x-x'|
$
\ \\

\noindent
$\quad$for all $x,x'\in\bar{x}(t)+\bar \delta B,\, v\in F(t,x),\,\mathrm{a.e.}\; t\in[S,T].$
\ \\

\noindent
{\bf (HM3):}
$e(.,.)$ is Lipschitz continuous on bounded subsets
\ \\

\noindent
Then 
there exist $p(.)\in W^{1,1}\left([S,T];\m R^n\right)$
and $\lambda \geq 0$ such that
\begin{itemize}
\item[(i)] $(p(.), \lambda)\; \not= \;(0,0)$,
\item[(ii)] $(-\dot p(t),\dot{\bar{x}}(t))\in \mbox{co}\,\partial_{x,p} H_{\lambda}^{i}(\bar{x}(t),p(t)) - \lambda\mbox{co}\, \tilde{\partial}_{x}e(t,\bar{x}(t))\times\{0\} \\
\\
 a.e. \; t \in [t_{i},t_{i+1}]$,\, $i=0,\ldots, N-1$,
where
$
\tilde{\partial}_{x}e(t,x)\,:=\, \underset{s \rightarrow t}{\mbox{lim sup}} \;\mbox{co}\, \partial_{x} e(s,x)\,,
$
\item[(iii)] 
$(p(S),-p(T)) \in  \lambda \partial g(\bar x(S),\bar x(T))+ N_{C}(\bar x(S),\bar x(T))$,
\end{itemize}
%

\noindent
Furthermore, there exists a piecewise Lipschitz continuous function $r(.):[S,T]\rightarrow \R$, with discontinuities confined to the set $\{t_{1}, \ldots, t_{N-1}\}$ satisfying
\ \\

\noindent
(iv) $\dot r(t)\;\in \; \lambda\mbox{co}\, \partial_{t} e(t,\bar{x}(t))\; \mbox{a.e.}$
\ \\

\noindent
(v) $ r(t)\,=\, H^{i}_{\lambda}(t, \bar{x}(t),p(t))\;  \mbox{a.e. } t \in [t_{i-1},t_{i}]$, $i=1,\ldots,N$.
\end{prop}

\noindent
{\bf Proof.} The proof is based on a well-known idea for deriving conditions on $H[.]$, when the data is Lipschitz continuous in time. It is to observe that the  given $L^{\infty}$-local minimizer $\bar x(.)$ remains an $L^{\infty}$-local minimizer after the domain of the optimization problem is enlarged to include the effects of an arbitrary change of independent variable, and then to apply standard necessary conditions to the modified optimization problem. We thereby  obtain, not just the usual necessary conditions for the $L^{\infty}$-local minimizer, but the desired additional  information about  $H[.]$. For problem $(P)'$, the data is not Lipschitz continuous in $t$. But  it is piecewise continuous, and we can use the idea to establish piecewise Lipschitz continuity of $H[.]$. Take an $L^{\infty}$ minimizer $\bar{x}(.)$ for $(P'')$ which is minimizint with respect to feasible $F$-trajectories that satisfy:
$$
||x(.)-\bar x(.)||_{L^{\infty}}\,\leq \, \bar{\epsilon}\,
$$
for some $\bar \epsilon >0$.
 Consider the 
function:
\begin{equation}
\label{solution}
\{(\bar z_{i}(.), \bar \tau_{i}(.), \bar x_{i}(.))\}_{i=0}^{N-1}: [S,T] \rightarrow (\R \times \R \times \R^{n})^{N} 
\end{equation}
in which, for $i=0,\ldots,N-1$,
$$
(\bar z_{i}(s), \bar \tau_{i}(s), \bar x_{i}(s))\;=\;
\left\{
\begin{array}{ll}
\left(\int_{t_{i}}^{s}L(s,\bar x(s'), \dot {\bar {x}}(s')) ds',  s, \bar x(s)\right)&
\mbox{for } s \in [t_{i},t_{i+1}]
\\
(0,t_{i},\bar x(t_{i}))& 
\mbox{for } s \in [S,t_{i})
\\
(\int_{t_{i}}^{t_{i+1}}L(s',\bar x(s'), \dot {\bar {x}}(s'))ds',t_{i+1},\bar x(t_{i+1}))& 
\mbox{for } s \in (t_{i},T]\,.
\end{array}
\right.
$$
We claim that, for $\epsilon' >0$ sufficiently small, (\ref{solution}) is an $L^{\infty}$-local minimizer for $(P'')$:
\noindent
\ \\

\noindent
$(P'')\:\left\{ \begin{array}{l}
\mathrm{Minimize }\, g(x_{0}(S),x_{N-1}(T))+ \sum_{i=0}^{N-1} z_{i}(T) 
\\
\mbox{over absolutely continuous functions }
\\
\hspace{0.3 in}\{( z_{i}(.),  \tau_{i}(.),  x_{i}(.))\}_{i=0}^{N-1}: [S,T] \rightarrow (\R \times \R \times \R^{n})^{N}, 
\\
 \mbox{ satisfying, for a.e. $t \in [t_{i},t_{i+1}] \mbox{ and } i =0,\dots,(N-1)$ },
\\ 
(\dot z_{i}(s), \dot \tau_{i} (s),\dot{x}_{i}(s)\,\in \, 
\\
\hspace{0.0 in}\left\{
\begin{array}{l}
\{(1+w)L_{i}(x_{i}(s),v), (1+w), (1+w)v\,|
\\
\hspace{0.5 in} w \in [1-\epsilon',1+\epsilon'], v \in F_{i}(s,x_{i}(s)) \} \;
\;
 \mbox{if } s \in [t_{i},t_{i+1}]
\\
(0,0,0)\hspace{2.4 in}\mbox{if } s \in [S,T]\backslash [t_{i},t_{i+1}]
\end{array}
\right.
\\
(x_{0}(S), x_{N-1}(T))\in C, \,
\\
(z_{i}(S), \tau_{i}(S), \tau_{i}(T))=(0,t_{i},t_{i+1} ) \mbox{ for }  i=0,\ldots, N-1
\\
x_{i}(T)=x_{i+1}(S) \mbox{ for } i=0,\ldots, N-2\,. 

\end{array}\right.
$
\ \\

\noindent
To see this, fix $\alpha >0$  and $\epsilon' \in (0,1)$, and take any $F$-trajectory for $(P'')$, 
$$
\{( z_{i}(.),  \tau_{i}(.),  x_{i}(.))\}_{i=0}^{N-1},
$$
 satisfying
$$
||\{( z_{i}(.),  \tau_{i}(.),  x_{i}(.))\} - \{\bar z_{i}(s), \bar \tau_{i}(s), \bar x_{i}(s))\}|| _{L^{\infty}} \,\leq \, \alpha\,.
$$
Then, in view of the dynamic constraint in $(P'')$, there exist measurable functions $w(.):[S,T] \rightarrow [1-\epsilon', 1+\epsilon'] $,  and $v(.):[S,T]\rightarrow \R^{n}$ such that
$$
\hspace{0.0 in}\left\{
\begin{array}{l}
\{(1+w)L(s, x_{i}(s),v), (1+w), (1+w)v\,|
\\
\hspace{0.5 in} w \in [1-\epsilon',1+\epsilon'], v \in F(s,x_{i}(s)) \} \;
\;
 \mbox{if } s \in [t_{i},t_{i+1}]
\\
(0,0,0)\hspace{2.4 in}\mbox{if } s \in [S,T]\backslash [t_{i},t_{i+1}]
\end{array}
\right.
\\
$$
for a.e. $s \in [t_{i},t_{i+1}] \mbox{ and } i =0,\dots,(N-1)$. Now consider the  function $\phi(.): [S,T] \rightarrow \R$ 
$$
\phi(s)\,=\, S+ \int_{[S, s]} \{\sum_{i=0}^{N-1} (1+w(s')) \chi_{[t_{i},t_{i+1}]}(s')\}ds'\;.
$$
In view of the constraints imposed in problem $(P'')$, $\phi(.)$ is a Lipschitz continuous, strictly increasing function with a Lipschitz continuous inverse, such that $\phi(t_{i})=t_{i}$ for  $i=0,\ldots,N$. Furthermore
$$
y(s) = \sum_{i=0}^{N-1} x_{i}(s)\chi_{[t_{i},t_{i+1})}(s)
$$
is an absolutely continuous function. It is a straightforward `change of independent variable' exercise to show that $x(.):[S,T] \rightarrow \R^{n}$ defined by 
$$
x(t) \,:=\, y( \phi^{-1}(t))
$$
is a feasible $F$-trajectory for $(P')$ and has the same cost as 
$\{( z_{i}(.),  \tau_{i}(.),  x_{i}(.))\}_{i=0}^{N-1}$ has for $(P'')$. Also, 
(\ref{solution}) has the same cost for $(P'')$. by choosing $\epsilon$ and $\alpha$ sufficiently small, we can arrange that, whatever the choice of $\{( z_{i}(.),  \tau_{i}(.),  x_{i}(.))\}$, we have $||x(.)-\bar x(.)||_{L^{\infty}} \leq \bar \epsilon$. The claim is confirmed. The assertions of proposition now follow from an application of known necessary conditions to $(P'')$ \cite[ThmXX]{vinter}.
\ \\

\noindent
We shall also require certain convergence properties, as the mesh size tends to zero, of interpolants of suitably bounded points  on a grid, summarized in the lemma:
\begin{lem}
\label{lem2.3}
Take  a compact set  $A \subset \R^{k}$, a continuous function $\bar{x}(.):[S,T]\rightarrow \R^{n}$ and a  multifunction $F(.,.,.): [S,T]\times \R^{n}\times A \rightarrow \R^{n}$ which has bounded variation along  $\bar{x}(.)$ uniformly over $A$.  Denote by $\eta(.)$  its cummulative variation function, and by $\eta_{\epsilon}^{\delta}(.)$ the related $(\delta, \epsilon)$-perturbed functions. 
Assume that hypotheses (C1) and (C2) of Prop. \ref{prop1} are satisfied.
Take a sequence of families of numbers $\{d_{j}^{i}\}_{j=1}^{N_{i}-1}$
and a sequence of numbers $\{m_{0}^{i}\}$.
For each $i$, let ${\cal T }_{i}=\{t_{0}^{i}, t_{1}^{i},\ldots,,t_{N_{i}-1}^{i}, t_{N_{i}}^{i}=T\}$ be a partition of $[S,T]$ and define the piecewise constant interpolation function
$$
m_{i}(t)\;=\;
m_{0}^{i}+ 
\sum_{j=1}^{N_{i}-2}d_{j}^{i}\chi_{[t_{j}^{i},t_{j+1}^{i})}(t)
+ 
d_{N_{i}-1}^{i}
\chi_{[t_{N_{i}-1}^{i},t_{N}^{i}] }(t)\;.
$$
Assume that
\begin{itemize}
\item[(A1):]
diam$\{ {\cal T}_{i} \} \,\rightarrow\, 0\;\mbox{ as }\; i \rightarrow \infty\,.$
\item[(A2):] $\{m_{0}^{i}\}$ is a bounded sequence.
\item[(A3):] There exist an integer valued function $(\epsilon,\delta)\rightarrow I(\epsilon, \delta)$ such that, for any $\epsilon >0$ and $\delta>0$,
$$
|d_{j}^{i}|\,\leq \, \eta^{\delta}_{\epsilon}(t^{i}_{j})- \eta^{\delta}_{\epsilon}(t_{j-1}^{i})\quad \mbox{ for all } j \in \{1,\ldots, N_{i}\},\; \mbox{and } i\geq I(\epsilon , \delta) \,.
$$
\end{itemize}
Then there exist a  normalized function of bounded variation $m(.): [S,T]\rightarrow\R^{n}$
 and a countable set ${\cal A} \subset (S,T)$,
 such that, along some subsequence,
$$
m_{i}(t) \rightarrow m(t)\quad \mbox{ for all } t \in [S,T]\backslash {\cal A} \,
$$
and
$$
|m(t)-m(s)| \,\leq\,  \eta^{*}(t)-\eta^{*}(s)  \mbox{ for all } [s,t] \subset [S,T]\,,
$$
in which $\eta^{*}(.)$ is the normalized cummulative variation 
$$
\eta^{*}(t)\,: =\,
\left\{
\begin{array}{ll}
\eta(t)& \mbox{if } t = S \mbox{ or  } T
\\
\eta(t^{+}) & \mbox{if } t \in (S,T)\,

\end{array}
\right.
$$
where $\eta(t^{+}):=\lim_{s\downarrow t} \eta(s)$\,.
\end{lem}

\noindent
A proof of the lemma appears in the Appendix.

\section{Proof of Theorem  \ref{thm:A} }

\noindent
Consider first the case  $L(.,.,.)\equiv 0$. Here, the Hamiltonian $H_{\lambda}(.,.,.)$, defined by (\ref {hamiltonian}), no longer depends on $\lambda$ and we write it simply $H(.,.,.)$.
\ \\

\noindent   By reducing the size of $\bar \delta>0$ in hypotheses $(H1)-(H3)$ we can arrange that $\eta_{F}^{\bar \delta}(T) < \infty$. We may  choose $\bar \epsilon> 0$ such that $\eta^{\bar \delta}_{F,\bar \epsilon}(T) < \infty$. Since $\bar{x}(.)$ is an $L^{\infty}$-local minimizer,
we may arrange (again by reducing the size of $\bar{\delta}>0$ if required) that $\bar{x}(.)$ minimizes $g(x(S),x(T))$ over all feasible $F$-trajectories for $(P)$ satisfying $||x(.)-\bar x(.)||_{L^{\infty}} \leq \bar \delta$.
\ \\

\noindent
We impose the following temporary hypothesis:
\ \\

\noindent
{\bf (C):} $t\rightarrow F(t,.)$ is right and left continuous at $S$ and $T$ respectively, in the following sense: for some $\delta \in (0, \bar \delta) $,
$$
\underset{s\downarrow S}{\lim} 
\underset{x \in \bar{x} (S)+\delta B}{\sup}
d_{H}(F(S,x), F(s,x))\,=\,0
\,, \;
\underset{t\uparrow T}{\lim} 
\underset{x \in \bar{x}(T)+\delta B}{\sup}
d_{H}(F(T,x), F(t,x))\,=\,0 \,.
$$ 
Take a sequence of positive integers $N_{i}\uparrow \infty $. Define $\epsilon_{i}=|T-S| / N_{i}$.  For each $i$, let  $\{t_{0}^{i}=S, t^{i}_{1}\ldots, t^{i}_{N_{i}-1}, t^{i}_{N_{i}}=T\}$ be the uniform partition of $[S,T]$ into $N_{i}$ subintervals. Define
$$ F_{i}(t,x)\;:=\;
\sum_{j=0}^{N_{i}-1}F(t_{j}^{i},x) \chi_{[t_{j}^{i},t_{j+1}^{i})}(t) + F(t_{N_{i}-1}^{i},x)
 \chi_{[t_{N_{i}-1}^{i},t_{N_{i}}^{i}]}(t) \;.
$$
In view of $(H2)$,
\begin{eqnarray*}
\int_{S}^{T}d_{F_{i}(t, \bar x(t))}( \dot{\bar x}(t))dt 
&\leq& \frac{T-S}{N_{i}}\times \left( \sum_{j=0}^{N_{i}-1} \eta^{\bar \delta}_{F, \bar \epsilon} (t^{i}_{j+1}) -\eta^{\bar \delta}_{F,\bar \epsilon}(t^{i}_{j})\right)
\\
&=& \frac{T-S}{N_{i}}\times \left( \eta^{\bar \delta}_{F,\bar \epsilon} (T) -\eta^{\bar \delta}_{F, \bar \epsilon}(S)\right),
\\
&\rightarrow& 0\,\; \mbox{as $i \rightarrow \infty\,$.}
\end{eqnarray*}
It follows from Filippov's Existence Theorem $[]$ that there exists an absolutely continuous function $z_{i}(.)$ such that $\dot z_{i}(t)\in F_{i}(t,z_{i}(t))$, a.e., and $z_{i}(S)= \bar{x}(S)$, and a sequence $\alpha_{i}\downarrow 0$ such that
$$
||z_{i}(.)- \bar{x}(.)||_{L^{\infty}}\,\leq\, \alpha_{i} \quad 
$$
for all $i$ sufficiently large. For each $i$ we have $h(z_{i}(t))\leq k_{h}\alpha_{i}$, where $k_{h}$ is the constant of $(H3)$ and $(z_{i}(S),z_{i}(T))\in \tilde C + \sqrt{2}\times \alpha_{i} \B$. Notice that
$$
\tilde C + \sqrt{2}\times \alpha_{i} \B \;\subset \; \{ (x_{0},x_{1})\,|\,  h(x_{0}) \vee h(x_{1}) \,\leq\,2 \times k_{h}\alpha_{i}\}\,.
$$
Take a sequence $\beta_{i} \downarrow 0$ such that  $\beta_{i}> 2k_{h}\alpha_{i}$ for each $i$ . Then
$$
\tilde C + \sqrt{2}\times \alpha_{i} \B \;\subset\; \{(x_{0},x_{1})\,|\,  h(x_{0})\vee h(x_{1})\,<\,\beta_{i}\}\,.
$$
Now take any sequence of numbers $K_{i}\uparrow \infty$ and, for each $i$, consider the optimization problem
\ \\

\noindent
$\:(P_{i})\left\{ \begin{array}{l}
\mathrm{Minimize }\, g(x(S),x(T)) + \int_{S}^{T}\left(|x(t)-\bar{x}(t)|^{2}+K_{i}\left(h(x(t))-\beta_{i} \right)^{+}\right)dt
\\
 \mbox{such that }
\\ 
\dot{x}(t)\,\in \,  F_{i}(t,x(t))\;\mathrm{a.e.,}
\\
(x(S), x(T))\in \tilde C + \sqrt{2} \times \alpha_{i} \B\;,
\\
|x(t)-\bar x(t)|\,\leq\, \bar{\delta}/2\quad \mbox{for all } t \in [S,T]\,.
\end{array}\right.
$
\ \\

\noindent
For each $i$,  $(P_{i})$ has a minimizer $x_{i}(.)$, since the data for this problem satisfy the standard conditions for existence of minimizers.  (The assumed convexity of the velocity sets $F(t,x)$ is crucial here.) 
Notice that the $\alpha_{i}$'s have been chosen to ensure existence of feasible $F$-trajectories for this problem. Since the $x_{i}(.)$'s are uniformly bounded and the $\dot x_{i}(.)$'s are uniformly integrably bounded we know that, along some subsequence (we do not relabel),
\begin{equation}
\label{4.1}
||x_{i}(.)-x'(.)||_{L^{\infty}} \,\rightarrow 0\,,
\end{equation} 
for some absolutely continuous function $x'(.)$ satisfying $||x'(.)-\bar x(.)||_{L^{\infty}}\,\leq\, \bar \delta/2$. Appealing once again to Filippov's Existence Theorem we can show that, for each $i$ sufficiently large, there exists an $F$-trajectory $y_{i}(.)$ such that
$||x_{i}(.)-y_{i}(.)||_{L^{\infty}} \rightarrow 0$. It follows from (\ref{4.1}) that
$$
||y_{i}(.)-x'(.)||_{L^{\infty}} \,\rightarrow \, 0\,.
$$
But then, by standard closure properties of solutions of convex valued differential inclusions, $x'(.)$ is an $F$-trajectory. Notice that, since $(x_{i}(S), x_{i}(T))\in \tilde C + \sqrt{2} \times \alpha_{i} \B$ for each $i$,
\begin{equation}
\label{4.2}
(x'(S),x'(T)) \,\in\, \tilde C\,.
\end{equation}
Observe next that the $z_{i}(.)$'s satisfy the conditions
$$
(z_{i}(S), z_{i}(T))\in \tilde C + \sqrt{2} \times \alpha_{i} \B,\, \underset{t \in [S,T]}{\mbox{max}} \left(h(z_{i}(t))-\beta_{i}\right)< 0\; \mbox {and}\;||z_{i}(.)-\bar{x}(.)||_{L^{\infty}}\leq \bar{\delta}/2\,,
$$ 
for $i$ sufficiently large. It follows that $z_{i}(.)$ is feasible for $(P_{i})$ and cannot have cost less than that of  $x_{i}(.)$. But then
\begin{eqnarray*}
&& g(x_{i}(S),x_{i}(T)) + \int_{S}^{T}|x_{i}(t)-\bar{x}(t)|^{2}dt +K_{i} \int_{S}^{T} \left( h(x_{i}(t))-\beta_{i} \right)^{+}dt
\\
&&\hspace{0.3 in}
\leq \; g(z_{i}(S),z_{i}(T)) +  \int_{S}^{T}|z_{i}(t)-\bar{x}(t)|^{2}dt + 0\,.
\end{eqnarray*}
Since  $||z_{i}(.)-\bar x(.)||_{L^{\infty}} \,\rightarrow 0$  and $||x_{i}(.)- x'(.)||_{L^{\infty}} \,\rightarrow 0$, we have
\begin{eqnarray}
\nonumber
&& g(x'(S),x'(T)) + \int_{S}^{T}|x'(t)-\bar{x}(t)|^{2}dt +
\underset{i \rightarrow \infty}{\mbox{lim sup}}\,
K_{i}\int_{S}^{T}\left( h(x_{i}(t))-\beta_{i} \right)^{+}dt
\\
&&
\label{4.3}
\hspace{3.5 in}
\leq \; g(\bar x(S), \bar x(T)) \,.
\end{eqnarray}
It follows from $K_{i}\uparrow \infty$ that 
$
\; \int_{S}^{T} \left( h(x_{i}(t))-\beta_{i} \right)^{+}dt\,\rightarrow 0\,\,.
$
\ \\

\noindent
We deduce from the continuity of $t \rightarrow h(x(t))$ that
$$
h(x'(t))\,\leq\, 0 \quad \mbox{for all } t \in [S,T]\,.
$$
But then,  by (\ref{4.2}),  $x'(.)$ is a feasible $F$-trajectory for the original problem $(P)$, which satisfies $||x'(.)-\bar x(.)||_{l^{\infty}}\leq \bar{\delta}/2$. In consequence then of the local optimality of $\bar{x}(.)$,
$$
g(x'(S),x'(T))\,\geq\, g(\bar x(S), \bar x(T))\,.
$$
We deduce from this relation and (\ref{4.3}) that $x'(.)=\bar x(.)$ Then, by (\ref{4.1}),
\begin{equation}
\label{4.4}
||x_{i}(.)-\bar{x}(.)||_{L^{\infty}}\,\rightarrow \, 0\,.
\end{equation}
and
$$
K_{i}\,\int_{S}^{T} \left( h(x_{i}(t))-\beta_{i} \right)^{+}dt \,\rightarrow\, 0\,.
$$
It follows from Egorov's Theorem that, after subtracting a subsequence, 
\begin{equation}
\label{4.4d}
K_{i}\, \left( h(x_{i}(t))-\beta_{i} \right)^{+} \,\rightarrow\, 0\quad \mbox{a.e.}
\end{equation}
In view of (\ref{4.4}), $x_{i}(.)$ is a $L^{\infty}$-local minimizer for  a modified version of $(P_{i})$, resulting from removal of the constraint `$||x(.)-\bar x(.)||_{L^{\infty}}\leq \bar \delta/2$'). The hypotheses are satisfied for the application of the `multistage' necessary conditions Prop. \ref{multistage} of the preceding section.  Write $q$ for the costate variable and $\lambda$ for the cost multiplier. The Hamiltonian for $(P_{i})$ is
$$
H_{i}(t,x,q)- \lambda |x-\bar{x}(t)|^{2}-\lambda K_{i}\left( h(x)- \beta_{i} \right)^{+}\,,
$$
in which 
$$
H_{i}(t,x,q)\,=\, \max\{q\cdot e\,|\, e \in F_{i}(t,x)  \}\,.
$$
The necessary conditions assert the existence, for each $i$, of $q_{i}(.) \in W^{1,1}([S,T];\R^{n})$ and a piecewise absolutely continuous function $r^{i}(.):[S,T]\rightarrow \R$, with possible jumps at $t_{1}^{i},\ldots, t^{i}_{N_{i}-1}$ and right continuous on $(S,T)$, such that
\begin{eqnarray}
&& (-\dot{q}_{i}(t),\dot{x}_{i}(t)) \in  \mbox{co}\, \partial_{x,p}H_{i}(t,x_{i}(t),q_{i}(t))
\label{4.4b}
\\
&&   
\hspace{1.5 in}
- 2 \lambda_{i}(x_{i}(t)-\bar x(t)) - 2 \gamma_{i}(t)K_{i}\nu_{i}(t) \quad \mbox{a.e.,}
\nonumber
\end{eqnarray}
in which $\gamma_{i}(.): [S,T]\rightarrow \R^{n}$ and $\nu_{i}(.): [S,T]\rightarrow [0,1]$ are measurable functions such that
\begin{equation}
\label{4.4c}
\gamma_{i}(t)\,\in\, \mbox{co}\,\partial^{>}_{x}h(x_{i}(t))
\end{equation}
for almost all $t\in \{ s\in[S,T]: \nu_{i}(s)>0  \}.$
To derive these relations, we have used the subdifferential calculus rule
$$
\partial[h(x)-\beta,0]^{+}\,\subset\,
\left\{
\begin{array}{ll}
\{\nu \partial h(x)\,|\, \nu \in [0,1] \} &\mbox{for } h(x)-\beta \geq 0
\\
\{0\}& \mbox{otherwise}
\;.
\end{array}
\right.
$$
Also:
\begin{equation}
\label{4.5}
\dot r_{i}(t)\,=\,2 \lambda_{i}(x_{i}(t)-\bar{x}(t))\cdot \dot{\bar{x}}(t)\quad \mbox{a.e.}
\end{equation}
and
\begin{equation}
\label{4.6}
-r_{i}(t)+H_{i}(t,x_{i}(t),q_{i}(t))- \lambda_{i}|x_{i}(t)-\bar{x}(t)|^{2}-\lambda_{i}K_{i}[h(x_{i}(t))-\beta_{i}]^{+}\;=\;0\; \mbox{ a.e.}
\end{equation}
Furthermore
\begin{equation}
\label{4.7}
q_{i}(t)\cdot \dot{x}_{i}(t)\,=\, \underset{e \in F_{i} (t, x_{i}(t))}{\max}q_{i}(t)\cdot e\;, \mbox{ a.e.}
\end{equation}
and
\begin{equation}
\label{4.8}
(q_{i}(S), -q_{i}(T))\,\in\, \lambda_{i}\partial g(x_{i}(S),x_{i}(T))+ N_{\tilde C + \sqrt{2}\alpha_{i}\B}(x_{i}(S),x_{i}(T))
\end{equation}
Now define
\begin{equation*}
\mu_{i}(dt)\,:= \, \lambda_{i}K_{i}\nu_{i}(t)dt
\end{equation*}
and
$$
p_{i}(t)\,:=\, q_{i}(t)- \int_{[S,t)}\gamma_{i}(s)\mu_{i}(ds)
$$
it follows from the last two relations that $(q_{i}(.),\lambda_{i}) \not= (0,0)$.
We can arrange then, by scaling the multipliers, that
\begin{equation}
\label{4.9a}
||p_{i}(.)||_{L^{\infty}}+\lambda_{i}+ ||\mu_{i}(.)||_{\mbox{TV}} \, =\, 1\,. 
\end{equation}
We note from (\ref{4.5}) that $r_{i}(.)$ has the representation
\begin{equation}
\label{4.10}
r_{i}(t)\,=\, \tilde r_{i}(t)- \int_{S}^{t}\left( 2 \lambda_{i} (x_{i}(s)-\bar{x}(s))\cdot \dot{\bar{x}} (s)  \right)ds\,,
\end{equation}
where $\tilde r_{i}(.)$ is the piecewise constant function with distributional derivative expressed in terms of the jumps
$$
\Delta^{i}_{j} \, = \,r_{i}(t^{i+}_{j})- r_{i}(t^{i-}_{j})\quad j=1,\ldots, N-1
$$
as
\begin{equation}
\label{4.11}
\dot{\tilde r}_{i} (t)\,=\, \sum_{j=1}^{N-1} \Delta_{j}^{i}\, \delta(t- t^{i}_{j})\,.
\end{equation}
Here, $\delta(.)$ is the Dirac delta function. We see from (\ref{4.6}) that the jumps in $\tilde r_{i}(.)$ are
$$
 \Delta_{j}^{i}\,=\,H_{i}(t^{i+}, x_{i}(t^{i}_{j}), q_{i}(t^{i}_{j}) ) -
H_{i}(t^{i-}, x_{i}(t^{i}_{j}), q_{i}(t^{i}_{j}) )\;,
$$
which, in view of $(H3)$, can be estimated by
\begin{equation}
\label{4.12}
| \Delta_{j}^{i}| \,\leq\, ||{q_{i}}(.)||_{\mbox{TV}}\times \left[ \eta^{\delta_{i}}_{F, \epsilon_{i}} (t^{i}_{j})- \eta^{\delta_{i}}_{F, \epsilon_{i}}(t^{i}_{j-1})  \right]\,.
\end{equation}
Since $\lambda_{i}K_{i}[h(x_{i}(.))- \beta_{i}] \rightarrow 0$ in $L^{1}$ and 
$t \rightarrow H_{i}(t,x_{i}(t),q_{i}(t))-\lambda_{i} |x_{i}(t)-\bar{x}(t)|^{2}$, $i=1,2,\ldots$, is 
a uniformly bounded sequence of functions, we can conclude from (\ref{4.5}) that $\{r_{i}(.)\}$ is bounded with respect to the $L^{1}$ norm. Note however that, by (\ref{4.11}) and (\ref{4.12}), $\tilde r_{i}(.)$ has total variation bounded by $||q_{i}(.) ||_{L^{\infty}} \times [\eta^{\delta_{i}}_{F,\epsilon_{i}}(T)- \eta^{\delta_{i}}_{F,\epsilon_{i}}(S)]$ A simple contradiction argument based on (\ref{4.10}) permits us to conclude that $r_{i}(S)(= \tilde r_{i}(S))$ is a bounded sequence.
\ \\

\noindent
We now apply Lemma \ref{lem2.3}, when we identify $m_{i}(.)$ with $\tilde r_{i}(.)$, $d^{i}_{j}$ with $r(t^{i+}_{j})-r(t^{i-}_{j})$ and $\eta^{\delta}_{F,\epsilon}(.)$ with $K' \eta^{\delta}_{F,\epsilon}(.)$, in which $K'$ is any number such that
$$
K' \, >\, \underset{i\rightarrow \infty}{\mbox{lim sup}}\; ||q_{i}(.)||_{L^{\infty}}\;.
$$
The lemma tells us that there exists a normalized function of bounded variation $r(.): [S,T]\rightarrow \R$ and a countable subset ${\cal A}\subset (S,T)$ such that 
$$
\tilde r_{i}(t) \rightarrow r(t)\; \mbox{for all } t \in [S,T]\backslash {\cal A}\,,
$$
and 
\begin{equation}
\label{4.16b}
|r(t)- r(s)|\,\leq\, K'  \times(\eta_{F}^{*}(t)-\eta_{F}^{*}(s))
\end{equation}
for all $[s,t] \subset [S,T]$.
\ \\

\noindent
We deduce from (\ref{4.6}), the facts that $h(x_{i}(S))<\beta_{i}$ and  $h(x_{i}(T))<\beta_{i}$, and the interim `continuity' hypothesis $(C)$ that, for some $\rho_{i} \downarrow 0$, 
\begin{eqnarray}
\nonumber
&&r_{i}(S)= H(S,x_{i}(S),p_{i}(S))+ \lambda_{i}|x_{i}(S)-\bar x(S)|^{2}\;
\\
\label{4.16c}
&&r_{i}(T)\in 
H(T,x_{i}(T),q_{i}(T))+\lambda_{i}|x_{i}(T)-\bar x(T)|^{2} + \rho_{i} \times
 ||q_{i}||_{L^{\infty}} \B\,.
\end{eqnarray}
(We have used here the fact that $F_{i}(S,.)=F(S,.)$, but  $F_{i}(T,.)=F(t,.)$ for some $t \in [S,T]$ such that $T-s\leq \epsilon_{i} =(T-S)/N^{i}$ for each $i$, since $F_{i}(.,.)$ is constructed from $F(.,.)$ by constant extrapolation from the left.)
The  $p_{i}(.)$'s are uniformly bounded and have a common Lipschitz constant.  The functions $\gamma_{i}(.) $  and $q_{i}(.)$ are uniformly bounded and also uniformly bounded in total variation, and the $\lambda_{i}$ are uniformly bounded. It follows that, for a subsequence,
\begin{eqnarray*}
&&p_{i}(.)\rightarrow p(.)\; \mbox{uniformly},\; \dot p_{i}(.) \rightarrow  \dot p(.) \mbox{ weakly in $L^{1}$ } \; \mu_{i}(.)\rightarrow \mu(.)\; \mbox{weakly*}
\\
&&  \gamma_{i}d\mu_{i} \rightarrow \gamma d \mu\; \mbox{weakly*} \;\mbox{and}\; q_{i}(.) \rightarrow q(.)  \; \mbox{weakly*} \;\;,
\end{eqnarray*}
for some Lipschitz continuous function $p(.)$, $\lambda \geq 0$ , some function of bounded variation $q(.)$, some measure $\mu(.)$. and some Borel measurable function $\gamma(.)$. Define
\ \\

\noindent
$$
q(t)\,:=\, p(t)+ \int_{[S,t)}\gamma(s)\mu(ds)
$$
A straightforward extension of the convergence analysis in  
\cite{vinter}
permits us pass to the limit in the relations  (\ref{4.4b}),  
(\ref{4.4c}), (\ref{4.8}) 
and (\ref{4.9a}) and thereby deduce
\begin{eqnarray}
\label{one1}
&&(-\dot{p}(t),\dot{\bar x}(t)) \in
 \,\tilde \partial_{x,p}H(t,\bar x(t),q(t)) \;\mbox{a.e.,}
\\
\label{two2}
&&
\lambda + ||p(.)||_{L^{\infty}} +  ||\mu(.)||_{\mbox{TV}}=1\,,
\\
\label{three3}
&&(p(S),-q(T)) \in  \lambda \partial g(\bar x(S),\bar x(T))+ N_{\tilde{C}}(\bar x(S),\bar x(T))\,,
\\
\label{four4}
&&
m(t)\in  \mbox{co}\, \partial_{x}^{>}h(t,\bar{x}(t))\quad \mu\mbox{-}\mathrm{a.e.}\; t\in[S,T]\,,
\\
\label{five5}
&&r(S)= H(S,\bar x(S),p(S))
\quad \mbox{and}\quad   r(T)= H(T,\bar{x}(T),q(T))\,.
\end{eqnarray}
\noindent
(Note that, for any $t \in [S,T]$, $H_{i}(t,.,.)$ and  $(H(t,.,. )$ may fail to coincide. Nonetheless, for any  index value $i$ and $t \in [S,T]$, $H_{i}(t,.,.)= H(s,.,. )$  for some $s$ such that $|t-s|\leq \delta_{i}$. So the partial subdifferential employed to capture limiting behavior for a particular time $t$, must take account of partial subgradients of $H(.,.,.)$  at neighbouring times $s$. This is why the  `hybrid' partial subdifferential  $\tilde \partial_{x,p} H$  appears in (\ref{one1}) in place of the customary $\mbox{co}\, \partial_{x,p} H$. To derive the relation (\ref{one1}) we make use of the fact: at all points $t$ in a subset of $[S,T]$ of full measure and some $\epsilon' >0$,  $t \rightarrow F(t,x)$ is continuous, uniformly over $x \in \bar{x}(t)+ \epsilon' \B$. See Prop. \ref{prop1}.)
\ \\

\noindent
We deduce from (\ref{4.6}), with the help of (\ref{4.4d}) that
$$
r(t)\,=\, H(t, \bar{x}(t), q(t)) \quad \mbox{a.e.}
$$
Reviewing the preceding relations (\ref{4.16b}) and (\ref{one}) - (\ref{five}) and, we see that the proof of the theorem in the `$L(.,.,.)\equiv 0$'  case is almost complete. But there are some minor matters that require attention, relating to the function $r(t)$. First, in the proof $K$ is taken to be any number  $K'>K$ where
$$
K\,=\, \underset{i \rightarrow \infty}{\mbox{lim sup}} 
||q_{i}(.) ||_{L^{\infty}}\,(\,=\, ||q(.) ||_{L^{\infty}} )
$$ 
To justify replacing $K'$ by $K$ in relation (\ref{4.16b}) we take a sequence $K_{j} \downarrow K$.  For each  $j$, we obtain the above relations with multipliers indexed by $j$. The desired necessary conditions, involving $K$, are obtained by extracting subsequences and passage to the limit in these relations.
\ \\

\noindent
The other matters concern the imposition of the temporary hypothesis $(C)$ and also the fact that the theorem statement additionally asserts the continuity of $r(.)$ at the endpoints of the interval $[S,T]$. Suppose condition $(C)$ is not valid. Then,  since $t\rightarrow F(t, .)$ is assumed to have bounded variation along $\bar{x}(.)$, $(C)$ is satisfied when the Hamiltonian inclusion $\dot x(t) \in  F(t,x(t))$ is replaced by  $\dot x(t) \in \tilde F(t,x(t))$, given by (\ref{modified}). Notice that $F(.,.)$ is obtained by changing $t \rightarrow F(t,.)$ only on a null-set,  and so this procedure does not affect $F$-trajectories. $\bar x(.)$ remains an $L^{\infty}$-local minimizer. The only respect in which this changes the preceding relations is to replace $\eta^{*}(.)$ in (\ref{4.16b}) by  the normalized cummulative variation function $\tilde  \eta^{*}(.)$ of $t \rightarrow \tilde F(t,.)$ in relation (\ref{4.16b}). We may  deduce from  
Prop. \ref{prop2.2}, however, that 
$$
r(t)-r(s) \,\leq\,\tilde \eta^{*}(t) - \tilde \eta^{*}(s)\,=\,  \eta^{*}(t) - \eta^{*}(s)
$$
for $[s, t] \subset (S,T)$.  But we also know from Prop. \ref{prop2.2} that $\tilde \eta^{*}(.)$ is continuous at the two endpoints of the interval $[S,T]$. So 
$$
\underset{t\downarrow S}{\lim}\,|r(t) -r(S)| \leq \underset{t\downarrow S}{\lim}\,\tilde \eta^{*}_{F}(t)- \eta^{*}_{F}(S^{+})=0
$$
and
$$
\underset{t\uparrow T}{\lim}\,|r(T) -r(t)| \leq \tilde \eta^{*}_{F}(T^{-})- \underset{t\uparrow T}{\lim}\,\eta^{*}_{F}(t)=0\;.
$$
We have shown that $r(.)$ has the desired continuity properties. This completes the proof of the special case of the theorem.
\ \\

\noindent
Now suppose that $L(.,.,.)$ is non-zero.  Then the assertions of the theorem may be deduced from those of the special case treated above by the well known state augmentation technique, based on the fact that
$(\bar{x}(.), \bar z(.))$ is an $L^{\infty}$-local minimizer for the optimal control problem, with state dimension $n+1$,
\ \\

\noindent
$(PA)\:\left\{ \begin{array}{l}
\mathrm{Minimize }\, g(x(S),x(T))+z(T)\\
\mbox{over absolutely continuous functions } (x(.),z(.)):[S,T]\rightarrow \mathbb{R}^{n+1} \mbox{ s.t. }
\\ 
(\dot{x}(t),\dot{z}(t) )\,\in \,  \tilde F(t,x(t))\;\mathrm{a.e.,}
\\
h(x(t))\,\leq\, 0 \quad \mbox{for all } t \in [S,T]
\\
(x(S), x(T))\in C\;\mbox{ and }\; z(S)=0\,.
\end{array}\right.
$
\ \\

\noindent
Here $\tilde F(.,.):\R^{n} \rightarrow \R^{n+1}$ is the multifunction:
$$
\tilde F(t,x)\;:=\; \left\{
 (e,\beta)\in \R^{n+1}\,|\, 
e \in F(t,x)\mbox{ and } L(t,x,e)\leq \beta \leq c' 
\right\}\,.
$$
in which $c'$ is any number satisfying $c' > c$, where $c$ is as in (H2). We  reproduce the preceding analysis in the proof of the  theorem, now with reference to $(\bar x(.), \bar z(.))$, interpreted as an $L^{\infty}$-local minimizer for 
$(PA)$.  This is permitted because all the relevant hypotheses are satisfied and since $(PA)$ has no integral cost term. 
The approximating `higher dimensional' costate arcs $(p_{i}(.), p^{0}_{i}(.))$  and associated arcs $(q_{i}(.), q_{i}^{0}(.))$,  which have an extra component, take the form $(p_{i}(.), p^{0}(.)) =(p_{i}(.),-\lambda_{i})$ and $ (q_{i}(.),-\lambda_{i})$. 
\ \\

\noindent
However we refine the analysis in one minor respect: this is to replace the estimate (\ref{4.12}) by the more refined relation,
$$
| \Delta_{j}^{i}| \,\leq\, ||{q_{i}}(.)||_{\mbox{TV}}\times \left[ \eta^{\delta_{i}}_{F, \epsilon_{i}} (t^{i}_{j})- \eta^{\delta^{i}}_{F,\epsilon_{i}}(t^{i}_{j-1})  \right] + \lambda_{i}\times \left[ \eta^{\delta_{i}}_{L, \epsilon_{i}} (t^{i}_{j})- \eta^{\delta^{i}}_{L,\epsilon_{i}}(t^{i}_{j-1})  \right]\,.
$$
in which the contributions to the estimates of the variation of the Hamiltonian on $[t^{i}_{j}, t^{i}_{j-1}]$), from the two components of 
$(q_{i}(.), (q_{i}^{0}(.)=\lambda_{i}))$, are now separated out. Here, $\eta^{\delta_{i}}_{F, \epsilon_{i}}$ and $\eta^{\delta_{i}}_{L, \epsilon_{i}}$ are the perturbed cummulative variation functions of $F(.,.)$ and $L(.,.,.)$ respectively.
\ \\

\noindent
We thereby arrive at a set of conditions from which may be deduced all the assertions of the theorem, when the integral cost term is present, with reference to  the $L^{\infty}$-local minimizer $\bar x(.)$ for $(P)$. 
\begin{center}
{\Large {\bf Appendix}}
\end{center}

\noindent
{\it Proof of Proposition \ref{prop1}:}\/
(a):
We prove the first assertion. Proof of the second assertion is analogous. Take any $\epsilon >0$ such that $\eta^{\bar \delta}_{\epsilon}(T) < + \infty$. Fix $\delta \in (0,\bar{\delta})$. Take any $x \in \bar{x}(\bar{s})+ \delta \B$, $a\in A$ and
$$
v \in \underset{s\downarrow \bar{s}}{\mbox{lim sup}}\, F(s,x,a)\;.
$$
By definition of the `lim sup', there exists $s_{i} \downarrow \bar{s}$ and $v_{i}\rightarrow  v$ such that
$$
v_{i}\in F(s_{i},x,a)    \mbox{ for all }i \;\mbox{ and }\;    v_{i}\rightarrow v \mbox{ as }  i\rightarrow \infty\,.
$$ 
The assertion (a) will follow if we can show that, also,
\begin{equation}
\label{lim inf}
v \in \underset{s\downarrow \bar{s}}{\mbox{lim inf}}\, F(s,x,a)\;,
\end{equation}
i.e. the `lim sup' and `lim inf' coincide, in which case the limit exists. To show (\ref{lim inf}) we  take an arbitrary sequence $t_{j}\downarrow \bar{s}$. By eliminating elements in the sequence $\{(s_{i},v_{i})\}$, we can arrange that, for every $j$,
$ \bar{s }\leq s_{j} < t_{j},\, t_{j}-\bar s\leq \epsilon \mbox{ and } x \in \bar{x}(t)+ \bar \delta \B\;\mbox{ for all } t \in [\bar{s},{t}_{j}]
$, $j =1,2,\ldots$, for all $a\in A$. But then, since $t_{j}-s_{j}\leq \epsilon$ and by an earlier listed `elementary' property of the $(\delta,\epsilon)$-perturbed cummulative variation function,
$$
d_{H}(F(t_{j},x,a), F(s_{j},x,a))\,\leq\, \eta_{\epsilon}^{\bar{\delta}}(t_{j})- \eta_{\epsilon}^{\bar{\delta}}(s_{j})\;.
$$
This means that, for each $j$, there exists $w_{j} \in F(t_{j},x,a)$ and 
$$
|v_{j}-w_{j}|\,\leq\, \eta_{\epsilon}^{\bar{\delta}}(t_{j})- \eta_{\epsilon}^{\bar{\delta}}(s_{j})\;.
$$
We know however that, since $\eta_{\epsilon}^{\bar \delta}(.)$ is a finite valued, monotone function, it has a right limit $\eta_{\epsilon}^{\bar \delta}(\bar{s}^{+})$ at $\bar{s}$. Hence
$$
\lim_{j\rightarrow \infty} |v_{j}-w_{j}|\,\leq\, 
\lim_{j \rightarrow \infty}\left( \eta_{\epsilon}^{\bar{\delta}}(t_{j})- \eta_{\epsilon}^{\bar{\delta}}(s_{j})\right) \,\leq\,
\eta_{\epsilon}^{\bar{\delta}}(\bar{s}^{+})- \eta_{\epsilon}^{\bar{\delta}}(\bar{s}^{+})=0\;.
$$
It follows that $v_{j}-w_{j}\rightarrow 0$. But then $v= \lim_{j}v_{j}= \lim_{j} w_{j}$. Since $t_{j}\downarrow \bar s$ was an arbitrary sequence, we conclude (\ref{lim inf}). We have confirmed (a).
\ \\

\noindent
(b) These assertions follow from (a), together with the compactness of the set $A$ and of the $\delta$ balls about $\bar{x}(\bar{s})$ and $\bar{x}(\bar{t})$, and with the assumed continuity properties of $(x,a)\rightarrow F(t,x,a)$.
\ \\

\noindent
(c) Let ${\cal A}$ be the empty or countable subset of $(S,T)$ comprising points at which the finite-valued monotone function $\eta_{\bar \epsilon}^{\bar \delta}(.)$ is discontinuous. Fix a point $t \in (S,T)\backslash {\cal A}$, $a\in A$ and $x \in \bar{x}(t)+\delta \B$. Take any $\rho >0$. Since $\eta_{\epsilon}^{\bar \delta}(.)$ is continuous at $t$, we may choose $\gamma >0$ such that 
$$
\eta_{\epsilon}^{\bar \delta}(t+\gamma)-\eta_{\epsilon}^{\bar \delta}(t-\gamma) \,\leq\, \rho\;.
$$
So, for any $t' \in [S,T]$ such that $|t'-t| \leq \rho \wedge \bar \epsilon$, 
\begin{eqnarray*}
\sup\{d_{H}(F(t',x,a),F(t,x,a))\,|\,x\in \bar{x}(t)+ \bar \delta \B,\, a\in A \}&\leq&  \eta^{\bar \delta}(t' \vee t)-\eta^{\bar \delta}(t' \wedge t)
\\
&\leq& \eta^{\bar \delta}(t+\gamma)-\eta^{\bar \delta}(t-\gamma)\leq\rho\;.
\end{eqnarray*}
The continuity properties of $F(.,x,a)$ at $t$ have been confirmed.
\ \\

\noindent
{\it Proof of Prop. \ref{prop2.2}:}\/ Let $\hat{F}(.,.): [S,T]\times \R^{n}\times A\rightarrow \R^{n}$ be a multifunction such that, for $(t,x,a)\in [S,T]\times \R^{n}\times A$,
$$
\hat{F}(t,x,a)\,=\,
\left\{
\begin{array}{ll}
F(S^{+},x,a)& \mbox{if } t=S, \,  x \in \bar x(S)+ \delta \B \, \mbox{ and } \, a\in A
\\
F(t,x,a)& \mbox{otherwise . }
\end{array}
\right.
$$
(Note that $\hat F(t,x,a)$ is a modified version of $ F(t,x,a)$ that differs only at the left endpoint $t=S$.) Write $\hat{\eta}^{\delta}_{\epsilon}(.)$ to denote the $(\delta,\epsilon)$-perturbed cummulative variation function of $t \rightarrow \hat F(t,.,.)$ along $\bar{x}(.)$, for $a\in A$.
\ \\

\noindent
Fix $t \in (S,T)$. Take $\delta \in (0, \bar \delta)$ and $\epsilon >0$ such that $\eta^{\delta}_{\epsilon}(T) < + \infty$.
Let ${\cal T}= \{t_{0}=S,\ldots, t_{N}=t\}$ be an arbitrary partition of $[S,T]$. Take an arbitrary sequence $s_{j}\downarrow S$. For $j$ sufficiently large, 
\begin{eqnarray*}
&&\eta^{\delta}_{\epsilon}(t)\,\geq\, \sup \left\{ d_{H}(F(s_{j},x,a) ,F(S,x,a) )\;|\; x \in \bar x([S,s_{j}])+\delta B, a\in A \right\}
\\
&&\hspace{1.0 in} + \sup \left\{ d_{H}(F(s_{j},x,a) ,F(t_{1},x,a) )\;|\; x \in \bar x([s_{j},t_{1}])+\delta B, \, a\in A \right\}
\\
&&\hspace{1.0 in} + \sum_{i=1}^{N-1}\sup \left\{ d_{H}(F(t_{i+1},x,a) ,F(t_{i},x,a) )\;|\; x \in \bar x([t_{i},t_{i+1}])+\delta B, a\in A \right\}\;.
\end{eqnarray*}
In view of Prop. \ref{prop1}, we may pass to the limit as $j\rightarrow \infty$ in this relation to obtain:
\begin{eqnarray*}
&&\eta^{\delta}_{\epsilon}(t)\,\geq\, \sup \left\{ d_{H}(F(S^{+},x,a) ,F(S,x,a) )\;|\; x \in \bar x(S)+\delta B,\, a\in A\right\}
\\
&&\hspace{1.0 in} + \sup \left\{ d_{H}(F(S^{+},x,a) ,F(t_{1},x,a) )\;|\; x \in \bar x([S,t_{1}])+\delta B, \, a\in A \right\}
\\
&&\hspace{1.0 in} + \sum_{i=1}^{N-1}\sup \left\{ d_{H}(F(t_{i+1},x,a) ,F(t_{i},x,a) )\;|\; x \in \bar x([t_{i},t_{i+1}])+\delta B, \, a\in A \right\}\;.
\end{eqnarray*}
Since the partition ${\cal T}$ was an arbitrary partition such that diam$\{ {\cal T}\} \leq \epsilon$ it follows that 
\begin{equation}
\label{one}
\eta^{\delta}_{\epsilon}(t) \geq \sup \left\{ d_{H}(F(S^{+},x,a) ,F(S,x,a) )\;|\; x \in \bar x(S)+\delta B, \, a\in A\right\} + \hat \eta^{\delta}_{\epsilon}(t)\,.
\end{equation}
Take again an arbitrary partition ${\cal T}= \{t_{0}=S, \ldots, t_{N}=t\}$ of $[S,t]$. Let $\bar \delta$, $\delta$ and $\epsilon$ be as before. We have
\begin{eqnarray}
\nonumber
&&\hat \eta^{\delta}_{\epsilon}(t)\,\geq\, \sup \left\{ d_{H}(F(S^{+},x,a) ,F(t_{1},x,a) )\;|\; x \in \bar x([S,t_{1}])+\delta B, a\in A \right\}
\\
\label{three}
&& + \sum_{i=1}^{N-1}\sup \left\{ d_{H}(F(t_{i+1},x,a) ,F(t_{i},x,a) )\;|\; x \in \bar x([t_{i},t_{i+1}])+\delta B, a\in A \right\}.
\end{eqnarray}
But by the triangle inequality we have, for each $x\in \{\bar x(t) +\delta \B\,|\,t \in [S,t_{1}]\} $ and $a\in A$,
$$
d_{H}(F(S^{+},x,a) ,F(t_{1},x,a) )\,\geq\, d_{H}(F(S,x,a) ,F(t_{1},x,a) ) - d_{H}(F(S^{+},x,a) ,F(S,x,a))\;.
$$
Furthermore,
\begin{eqnarray*}
&& \max d_{H}(F(S^{+},x,a) ,F(t_{1},x,a) )\geq
\\
&& \hspace{0.5 in}  \max d_{H}(F(S,x,a) ,F(t_{1},x,a) ) - \max\{ d_{H}(F(S^{+},x,a) ,F(S,x,a) )\,,
\end{eqnarray*}
where, in each term, the max is taken over $(x,a)\in \{\bar x(t) +\delta \B\,|\,t \in [S,t_{1}]\} \times A$.
Noting that ${\cal T}$ was an arbitrary partition such that diam$\{ {\cal T} \}\leq \epsilon$,  we deduce from (\ref{three})
that
\begin{eqnarray}
\nonumber
&&\hat \eta^{\delta}_{\epsilon}(t)\,\geq\, \eta^{\delta}_{\epsilon}(t) - \max\{ d_{H}(F(S^{+},x,a) ,F(S,x,a) )\,|\,(x,a)\in \{\bar x(t) +\delta \B\,|\,t \in [S,t_{1}]\}\times A    \}\,.
\end{eqnarray}
This relation combines with (\ref{one}) to yield
\begin{eqnarray*}
0\,\leq\,
\eta^{\delta}_{\epsilon}(t)-\hat \eta^{\delta}_{\epsilon}(t) - \max\{ d_{H}(F(S^{+},x,a) ,F(S,x,a) )\,|\,x\in \bar x(S)+ \delta \B, \, a\in A\}
\leq \Delta(\epsilon, \delta)\,,
\end{eqnarray*}
in which
\begin{eqnarray*}
&&  \Delta(\epsilon, \delta)\,:=\,
\max\{ d_{H}(F(S^{+},x,a) ,F(S,x,a) )\,|\,x\in \{\bar x(t) +\delta \B\,|\,  t \in [S,(S+\epsilon)\wedge T]\}, \, a\in A \}  
\\ 
&&
\hspace{1.8 in}
-\,\max\{ d_{H}(F(S^{+},x,a) ,F(S,x,a) )\,|\,x\in \bar x(S)+ \delta \B, a\in A\}\;.
\end{eqnarray*}
Since, as is easily shown, $F(S^{+},.,.)$ has modulus of continuity $\theta(.)$ on $\{\bar{x}(S)+\delta \B\} \times A$, where $\theta(.)$ is as in hypothesis (C2) and $\bar x(.)$ is continuous, 
$$
\lim_{\epsilon' \downarrow 0} \Delta(\epsilon', \delta)=0\;.
$$
We deduce that 
$$
\eta^{\delta}(t)=\hat \eta^{\delta}(t) + \max \{ d_{H}(F(S^{+},x,a) ,F(S,x,a) )\,|\,x\in \bar x(S) + \delta \B, \, a\in A\}\;.
$$
In the limit as $\delta \downarrow 0$ we obtain
\begin{equation}
\label{1.1}
\eta(t)=\hat \eta(t) +\underset{a\in A}{\sup} \,  d_{H}(F(S^{+},\bar{x}(S),a) ,F(S,\bar{x}(S),a) )\;.
\end{equation}
Since $\hat{\eta}(t)$ and $\tilde{\eta}(t)$ coincide for $t < T$, we have shown that
\begin{equation}
\label{1.1a}
\eta(t)=\tilde \eta(t) +\underset{a\in A}{\sup} \, d_{H}(F(S^{+},\bar{x}(S),a) ,F(S,\bar{x}(S),a) )\; \mbox{for }t \in [S,T)\,.
\end{equation}
A similar analysis to that above yields
\begin{equation}
\label{1.2}
\hat{\eta}(T) = \tilde \eta(T) +\underset{a\in A}{\sup} \, d_{H}(F(T^{-},\bar{x}(T),a) ,F(T,\bar{x}(T),a) )\;.
\end{equation}
Combining (\ref{1.1}) (in the case $t=T$) and (\ref{1.2}) yields
\begin{eqnarray}
\nonumber
&&{\eta}(T) = \tilde \eta(T) +\underset{a\in A}{\sup} \, d_{H}(F(S^{+},\bar{x}(S),a) ,F(S,\bar{x}(S),a) )
\\
\label{1.2b}
&& \hspace{1.5 in}
+\underset{a\in A}{\sup} \, d_{H}(F(T^{-},\bar{x}(T),a) ,F(T,\bar{x}(T),a) )\,.
\end{eqnarray}
It follows from (\ref{1.1a}) that 
\begin{equation}
\label{1.2a}
\tilde \eta(t) -\tilde \eta(s) =  \eta(t) - \eta(s)\;, \mbox{ for all } [s,t] \subset (S,T)\;.
\end{equation}
The remaining assertions of the proposition will follow from (\ref{1.1a}) and (\ref{1.2b}), if we can verify the two assertions: $\tilde \eta(.)$ is right continuous at $S$ and left continuous at $T$. We proof the first assertion. The proof of the second is similar.
\ \\

\noindent
 Suppose the first assertion is not true. Then there exists $\alpha > 0$ such that
$
\tilde \eta(t)-(\tilde \eta(S)=0) \geq \alpha \mbox{ for all $t \in [S,T]$}
$. Choose any $\delta \in (0,\bar \delta)$ and $\epsilon >0$ such that $\tilde \eta^{\delta}_{\epsilon}(t) < \infty$. Then
$$
\tilde \eta^{\delta}_{\epsilon}(t) \,\geq \, \alpha \; \mbox{for all } t \in [S,T]\;.
$$
Notice that the choice of  $\epsilon$ does not depend on the choice of $\alpha$. We can then impose that 
$$\Delta(\epsilon,\delta)< \alpha/8.$$ 
By Prop \ref{prop1}, we can find $\bar s >0 $ such that
$$
\max_{\underset{a\in A} {x \in \bar{x}(S)+ \delta B,}}d_{H}(F(\bar{s},x,a), F(S^{+},x,a)) \,\leq \, \alpha /4\,.
$$
By the properties of the supremum, we can choose a partition $\{s_{0}, ..., s_{N}\}$  of $[S, \bar s]$, of diameter at most $\epsilon$,  such that
\begin{eqnarray}
\nonumber
\tilde \eta^{\delta}_{\epsilon}(\bar{s}) &\leq&  \max_{\underset{a\in A}{x \in \bar{x}([S,s_{1}])+ \delta B} }d_{H}(F(s_{1},x ,a),F(S^{+},x,a ) ) + \alpha /4
+\Sigma_{2}
\\
&=& \Delta(\epsilon,\delta)+\max_{\underset{a\in A} {x \in \bar{x}(S)+ \delta B,}}d_{H}(F(\bar{s},x,a), F(S^{+},x,a)) + \alpha /4+\Sigma_{2}
\\
\label{inequality}
&\leq& \alpha/8 +  \alpha/4 + \Sigma_{2}+\alpha/4 \,=\, \Sigma_{2} + 5\alpha/8\,, 
\end{eqnarray}
where 
$$
\Sigma_{2}\,:=\, \sum_{i=1}^{N-1}\max_{\underset{a\in A}{x \in \bar{x}([s_{i},s_{i+1}])+ \delta \B,}} d_{H}(F(s_{i+1},x,a ),F(s_{i},x,a ) ) \,.
$$
But we can also choose a partition $\{t_{0}, \ldots,t_{M}\}$ of $[S,s_{1}]$ (which will have diameter not greater than $\epsilon$) such that
$$
\alpha \leq \tilde \eta^{\delta}_{\epsilon}(s_{1}) \leq \Sigma_{1}+ \alpha/4 \,,
$$
where
$$
\Sigma_{1}\,:=\, \sum_{i=0}^{N-1}\max_{\underset{a\in A}{x \in \bar{x}([t_{i},t_{i+1}])+ \delta \B}} d_{H}(F(t_{i+1},x ,a),F(t_{i},x,a ) ) \,.
$$
It follows that 
$$
\Sigma_{1} \geq 3\alpha/4\,.
$$ 
But since the concatenation of $\{t_{0}\ldots,t_{M}\}$ and $\{s_{1}\ldots,s_{N}\}$ is a partition of   $[S,\bar{s}]$, of diameter no greater than $\epsilon$, we know from the preceding inequality that
$$
\tilde{\eta}^{\delta}_{\epsilon}(\bar{s}) \geq \Sigma_{1} + \Sigma_{2} \geq  \Sigma_{2} + 3\alpha/4\,.
$$
But this contradicts (\ref{inequality}). The assertion has been confirmed.
\ \\

\noindent
{\it Proof of Lemma \ref{lem2.3}:}\/  Notice that 
\begin{equation}
\label{C2.3}
m_{i}(S)=m_{0}^{i} \;\mbox{ and }\;
||m_{i}(.)||_{TV}\;\leq\; \eta(T)-\eta(S)\; \mbox{ for all } i\;.
\end{equation} 
These relations follow from the definition of the $m_{i}(.)$'s and $(H3)$. 
\ \\

\noindent
Take sequences $\delta_{m} \downarrow 0 $ and $\epsilon_{n} \downarrow 0$.
For each $(m,n)$ let ${\cal B}_{m,n}$ be the (possibly empty,) countable set comprising points of discontinuity of the monotone function $\eta^{\delta_{m}}_{\epsilon_{n}}(.)$. Fix $(m,n)$.
\ \\

\noindent
For each $i$, define the function $\tilde m(.):[S,T]\rightarrow \R^{n}$, which can be interpreted as an interpolant of the values of $m(.)$ at grid points, as follows:
\begin{equation}
\label{C3.1}
\tilde m_{i}(t)\,:=\,m_{i}(t_{k}^{i})+ \left(m_{i}(t_{k+1}^{i})- m_{i}(t_{k}^{i})\right)\times 
\frac{\eta
^{\delta_{m}}_{\epsilon_{n}}(t)- \eta^{\delta_{m}}_{\epsilon_{n}}(t^{i}_{k})}
{\eta^{\delta_{m}}_{\epsilon_{n}}((t^{i}_{k+1}))- \eta^{\delta_{m}}_{\epsilon_{n}}(t^{i}_{k})}
\end{equation}
if $t \in [t_{k}^{i},t_{k+1}^{i})$, for some $k= 0,\ldots, N_{i}-1$. Set $\tilde m_{i}(T):= m_{i}(T)$.
\ \\

\noindent
(The rightside of (\ref{C3.1}) is interpreted as  $m_{i}(t_{k}^{i})$  if $\eta^{\delta_{m}}_{\epsilon_{n}}(t^{i}_{j+1})-\eta^{\delta_{m}}_{\epsilon_{n}}(t^{i}_{j})=0$.)
\ \\

\noindent
{\it Claim:} 
\begin{eqnarray}
&&
\label{C3.2}
\tilde m_{i}(S)=m_{i}(S) \quad \mbox{and}\quad \tilde m_{i}(T)=m_{i}(T)\; \mbox{for all } i
\\
&&
\label{C3.3}
\tilde m_{i}(t)-m_{i}(t) \,\rightarrow\,0 \quad \mbox{for all}\quad t \in [S, T]\backslash {\cal B}_{m,n} 
\\
&&
\label{C3.4}
|\tilde m_{i}(t)-\tilde m_{i}(s)| \,\leq\,\eta^{\delta_{m}}_{\epsilon_{n}}(t)-\eta^{\delta_{m}}_{\epsilon_{n}}(s)\quad 
\\
&& \hspace{1.0 in} \mbox{for all}\quad [s,t]\subset [S,T] \mbox{ and } i\geq I(\delta_{m}, \epsilon_{n}) \,.
\nonumber
\end{eqnarray}
We verify the claim. (\ref{C3.2}) follows from the fact that $m_{i}(.)$ and $\tilde{m}_{i}(.)$ coincide at mesh  points, which include $S$ and $T$. Consider next  (\ref{C3.3}).  Take any $t \in (S,T)\backslash {\cal B}$. 
Then
$$
|\tilde m_{i}(t)- m_{i}(t)| = |\tilde m_{i}(t)- \tilde m_{i}(t^{i}_{j_{i}})| \,\leq\,  \eta^{\delta_{m}}_{\epsilon_{n}}(t^{i}_{j_{i}+1})-\eta^{\delta_{m}}_{\epsilon_{n}}(t^{i}_{{j}_{i}}) \,,
$$ 
for each $i$, where $j_{i}$ is the unique index value satisfying $t \in [t^{i}_{j_{i}}, t^{i}_{j_{i}+1}  )$, by (H3). 
But, by (H1),
$t^{i}_{j_{i}} \rightarrow t$ and $t^{i}_{j_{i}+1} \rightarrow t$. So (\ref{C3.3}) follows from the fact that $t$ is a continuity point of $\eta(.)$. 
\ \\

\noindent
Consider finally assertion (\ref{C3.4}). Take any $[s,t] \subset [S,T]$. Then either
\begin{itemize}
\item[Case 1:]$[s,t]\subset [t^{i}_{j},t^{i}_{j+1}]] $ for some $j \in \{0,\ldots ,N_{i}-1\}$
\item[Case 2:]$s \in  [t^{i}_{j-1},t^{i}_{j}]$ and $s \in  [t^{i}_{k},t^{i}_{k+1}]$ for some $j,k$ such that  $j\leq k$.
\end{itemize}
We consider only Case 2. (Verifying Case 1 is similar, but simpler.) We have
\begin{equation}
\label{C3.6}
\tilde m_{i}(t)-\tilde m_{i}(s)\,=\,\tilde m_{i}(t^{i}_{j})-\tilde m_{i}(s)+ \sum_{l=j}^{k-1}
\left(m(t^{i}_{l+1}) - m(t^{i}_{l}) \right)+ \tilde m(t) - m(t^{i}_{k})
\end{equation}
(We have used here the fact that $\tilde m_{i}(.)$ and $m(.)$ coincide at points $t^{i}_{j}$, $j =1,\ldots,N_{i}$.)
But from the definition of $\tilde m_{i}(.)$, and in view of (H3),
\begin{eqnarray*}
&&|\tilde m(t^{i}_{j}) - \tilde m(s)|= (\eta^{\delta_{m}}_{\epsilon_{n}}(t^{i}_{j})-\eta^{\delta_{m}}_{\epsilon_{n}}(s))\times \frac{m(t^{i}_{j})-m(t^{i}_{j-1})}{\eta^{\delta_{m}}_{\epsilon_{n}}(t^{i}_{j})-\eta^{\delta_{m}}_{\epsilon_{n}}(t^{i}_{j-1})} \leq \eta^{\delta_{m}}_{\epsilon_{n}}(t^{i}_{j}) - \eta(s)\;.
\\
&&|\tilde m(t) - \tilde m(t^{i}_{k})|= (\eta^{\delta_{m}}_{\epsilon_{n}}(t)-\eta^{\delta_{m}}_{\epsilon_{n}}(t^{i}_{k}))\times \frac{m(t^{i}_{k+1})-m(t^{i}_{k})}{\eta^{\delta_{m}}_{\epsilon_{n}}(t^{i}_{k+1})-\eta^{\delta_{m}}_{\epsilon_{n}}(t^{i}_{k})} \leq \eta^{\delta_{m}}_{\epsilon_{n}}(t) - \eta^{\delta_{m}}_{\epsilon_{n}}(t^{i}_{k})\;.
\end{eqnarray*}
It follows from (\ref{C3.4}) and (H3) that
\begin{eqnarray*}
|\tilde m_{i}(t)-\tilde m_{i}(s)| &\leq& \eta^{\delta_{m}}_{\epsilon_{n}}(t^{i}_{j})-\eta^{\delta_{m}}_{\epsilon_{n}}(s)+\sum_{l=j}^{k-1}\eta^{\delta_{m}}_{\epsilon_{n}}(t^{i}_{l+1})- \eta^{\delta_{m}}_{\epsilon_{n}}(t^{i}_{l})+ \eta(t)-\eta(t^{i}_{k})
\\
&=& \eta^{\delta_{m}}_{\epsilon_{n}}(t)+ 0 + \ldots +0 -\eta^{\delta_{m}}_{\epsilon_{n}}(s)= \eta^{\delta_{m}}_{\epsilon_{n}}(t)-\eta^{\delta_{m}}_{\epsilon_{n}}(s)\,.
\end{eqnarray*}
We have confirmed (\ref{C3.4}) (in the Case 2) and thereby verified the claim.
\ \\

\noindent
In view of (H1), we can deduce from (\ref{C3.2}) and (\ref{C3.4}) that the total variation of elements in the sequence $\{\tilde m_{i}(.) \}$ are uniformly bounded and that their initial values are uniformly bounded. 
It follows that  there exists a normalized
function of bounded variation $m(.): [S,T]\rightarrow \R^{n}$
and a countable set $\tilde{\cal A}\subset (S,T)$ such that, for some subsequence,
\begin{equation}
\label{C5.1}
\tilde m_{i}(t) \rightarrow m(t)\quad \mbox{for all} \quad t \in [S,T]\backslash \tilde {\cal A}  \,.
\end{equation}
But then, by (\ref{C3.3}),
$$
m_{i}(t)\rightarrow m(t)\quad \mbox{ for all} \quad t \in [S,T] \backslash  \left( 
\tilde{{\cal A}}\cup {\cal B}_{m,n}  \right)\,,
$$
Up to this point the index values $m,n$ have been fixed. We now let them vary. Define the countable set
$$
{\cal A}\,:=\, \tilde {\cal A}\cup_{m,n}{\cal B}_{m,n}\;. 
$$
Passing to the limit as $i \rightarrow \infty$ in (\ref{C3.4}) we deduce that, for every 
$[S, T]\backslash {\cal A}$, we have
$$
|m(t)-m(s)|\leq \eta^{\delta_{m}}_{\epsilon_{n}}(t)-\eta^{\delta_{m}}_{\epsilon_{n}}(s)
$$
for all $[s,t]\subset [S,T]$ such that $s,t\notin (S,T)\backslash {\cal A}$. Fixing $m$ and passing to the limit as $n \rightarrow \infty$ yields
$$
|m(t)-m(s)|\leq \eta^{\delta_{m}}(t)-\eta^{\delta_{m}}(s)\,.
$$
Then, passing to the limit as  $m \rightarrow \infty$ in this relation, we obtain
$$
|m(t)-m(s)|\leq \eta(t)-\eta(s)\,.
$$
This inequality is valid, we recall, for all $[s,t] \in [S,T]$ such that $s,t \in [S,T]\backslash {\cal A}$.  It implies however
$$
|m(t)-m(s)|\leq \eta^{*}(t)-\eta^{*}(s)\,.
$$
for all subintervals $[s,t] \subset [S,T]$, since the regularized cummulative variation $\eta^{*}(.)$ coincides with $\eta(.)$ on the complement of a countable subset of $[S,T]$ that incluldes $\{S\} \cup \{T\}$ and since $m(.)$ and $\eta^{*}(.)$ are right continuous on $(S,T)$.
\ \\

\noindent
{\it Proof of Lemma \ref{simple}:}\/ Write $\eta^{\delta}_{\epsilon, A}(.)$ for the $(\delta, \epsilon)$ perturbed cummulative variation with respect to $A$, etc. TTake $\bar{\delta}>0$ and $\bar{\epsilon} >0$ such that $\eta^{\delta}_{\epsilon, A}(.) < \infty$ for every $\delta \in (0,\bar{\delta})$, $\epsilon \in (0, \bar{\epsilon}]$. Take also any $[s,t] \subset [S,T]$. We shall show that
\begin{equation}
\label{2star}
\eta^{\delta}_{A_{1}}(t)-\eta^{\delta}_{A_{1}}(s)\;\leq\; \eta^{\delta}_{A}(t)-\eta^{\delta}_{A}(s)\;.
\end{equation}
Since this relation is valid for all $\delta>0 $ sufficiently small, we can deduce (\ref{1star}) by passing to the limit as $\delta \downarrow 0$.
\ \\

\noindent
If $s=S$ then (\ref{2star}) is obvisions, since $\eta^{\delta}_{A_{1}}(S)=\eta^{\delta}_{A}(S)=0$ and since $\eta^{\delta}_{A_{1}}(t)\leq \eta^{\delta}_{A}(t)$. The latter inequality is true because the left side involves taking the supremum over a small set, as compared with the right side. So we may assume that $s > S$.
\ \\

\noindent
Define $\eta^{\delta}_{\epsilon, A}([s,t])$ to be the $(\delta, \epsilon)$ perturbed cummulative variation with respect to $A$, when the underlying time interval is changed from $[S,T]$ to the smaller set $[s,T]$. It is clear that
\begin{equation}
\label{3star}
\eta^{\delta}_{\epsilon, A_{1}}(t)\;\leq\; \eta^{\delta}_{\epsilon,A}(t)\;.
\end{equation}
because the left side involves taking the supremum over a small set, in relation to the right side.
\ \\

\noindent
By consideration of arbitrary partitions of $[S,s]$ of diameter not greater than $\epsilon$, as well as their extensions to form partitions of the larger interval $[S,t]$, we deduce that
\begin{equation*}
\eta^{\delta}_{\epsilon,A}(t)-\eta^{\delta}_{\epsilon, A}(s)\;\geq\; \eta^{\delta}_{\epsilon, A}([s,t])\;.
\end{equation*}
Passing to the limit as $\delta \downarrow 0$ yields
\begin{equation}
\label{4star}
\eta^{\delta}_{A}(t)-\eta^{\delta}_{A}(s)\;\geq\; \eta^{\delta}_{A}([s,t])\;.
\end{equation}
On the other hand,  consideration of arbitrary partitions of $[S,t]$ of diameter not greater than $\epsilon$, as their refinements to include the intermediate point $s$,  we deduce that there exists  continuity modulus $\gamma(.)$ such that
\begin{equation*}
\eta^{\delta}_{\epsilon,A}(t)-\underset{t' \uparrow s}{\lim}\,\eta^{\delta}_{\epsilon, A}(t')\;\leq\; \eta^{\delta}_{\epsilon,A}([s,t]) + \gamma(\epsilon)\;.
\end{equation*}
Noting that  $\lim_{t' \uparrow s}\,\eta^{\delta}_{\epsilon, A}(t')\leq \eta^{\delta}_{\epsilon, A}(s) $  and passing to the limit as $\delta \downarrow 0$ yields
\begin{equation}
\label{5star}
\eta^{\delta}_{A}(t)-\eta^{\delta}_{ A}(s)\;\leq\; \eta^{\delta}_{A}([s,t])\;.
\end{equation}
(\ref{2star}) now follow from (\ref{3star}), (\ref{4star}) and (\ref{5star}).

\end{document}